\documentclass[11pt]{article}

\usepackage{amsmath,amssymb,amsfonts,amsthm}
\usepackage{moreverb,rotating,graphics}
\usepackage[T1]{fontenc}
\usepackage[latin1]{inputenc}
\usepackage{epsfig}
\usepackage[francais,english]{babel} 
\usepackage{color}
\usepackage{stmaryrd}
\usepackage{dsfont}
\usepackage{float}
\newtheorem{theorem}{Theorem}
\newtheorem{Prop}[theorem]{Proposition}
\newtheorem{remark}{Remark}
\newtheorem{cor}[theorem]{Corollary}
\newtheorem{defi}[theorem]{Definition}
\newtheorem{lemma}[theorem]{Lemma}

\def\NN{\mathbb{N}}
\def\ZZ{\mathbb{Z}}
\def\RR{\mathbb{R}}

\numberwithin{equation}{section}

\title{Elements of mathematical foundations for a numerical approach for weakly random homogenization problems} 
\author{A. Anantharaman$^{1}$ \thanks{The author acknowledges EADS IW for financial support.} \, and C. Le Bris$^{2}$\\ \\
       \footnotesize{Universit\'e Paris-Est, CERMICS, Project-team
         Micmac, INRIA-Ecole des Ponts,} \\  
\footnotesize{6 \& 8 avenue Blaise Pascal,
         77455 Marne-la-Vall\'ee Cedex 2, France} \\ 
\footnotesize{$^{1}$ \tt ananthaa@cermics.enpc.fr, \; $^{2}$ lebris@cermics.enpc.fr}\\
}

\begin{document}

\maketitle

\begin{abstract}
\noindent
This work is a follow-up to our previous work~\cite{ALB1}. It extends
and complements, both theoretically and experimentally, the results
presented there. Under consideration is the  homogenization of a model of a
weakly random heterogeneous material. The material consists of a reference periodic material
 randomly perturbed by another periodic material, so that its homogenized behavior is close to
 that of the reference material.  
 We consider laws for the random perturbations more
 general than in \cite{ALB1}. We  prove the validity of an asymptotic
 expansion in a certain class of settings. We also extend the formal approach introduced in \cite{ALB1}.
Our perturbative
approach shares
common features with a defect-type theory of solid state physics. The
computational efficiency of the approach is demonstrated.
\end{abstract}

{\footnotesize{{\bf Keywords}: Homogenization; Random Media; Defects}\\

}

\section{Introduction}

Our purpose is to  follow up on our previous study~\cite{ALB1}. Let us
recall, for consistency, that we consider homogenization for the following elliptic problem
\begin{eqnarray} \label{problemintro}
\left \{
\begin{aligned}
&-\mathrm{div} \left((A_{per}(\frac{x}{\epsilon}) + b_{\eta}(\frac{x}{\epsilon}, \omega) C_{per}(\frac{x}{\epsilon})) \nabla u_{\epsilon}\right) = f(x) \; \mathrm{in} \; \mathcal{D} \subset{\mathbb{R}^d}, \\
&u_{\epsilon} = 0 \; \mathrm{on} \; \partial \mathcal{D},
\end{aligned}
\right.
\end{eqnarray}
where the tensor $A_{per}$ models a reference $\mathbb{Z}^d$-periodic
material which is randomly perturbed by the $\mathbb{Z}^d$-periodic
tensor $C_{per}$, the stochastic nature of the problem being encoded in the
stationary ergodic scalar field $b_{\eta}$ (the latter getting small when  $\eta$ vanishes). 
We have studied in \cite{ALB1} the case of a perturbation that has a Bernoulli law with parameter $\eta$,
meaning that $b_{\eta}$ is equal to $1$ with probability $\eta$ and $0$ with
probability $1-\eta$. In the present work, we address more general
laws. The common setting is that all the perturbations we consider are,
to some extent, rare events which, although rare, modify the homogenized
properties of the material. Our approach is a perturbative approach, and
consists in approximating the stochastic homogenization problem for
\begin{eqnarray*} 
A_{\eta}(x,\omega) = A_{per}(x) + b_{\eta}\left(x, \omega\right)C_{per}
\end{eqnarray*}
using the periodic homogenization problem for $A_{per}$. In
short, let us say that our main contribution is to derive an expansion
\begin{eqnarray} \label{expansintro}
A_{\eta}^* = A_{per}^* + \eta \bar{A_1^*} + \eta^2 \bar{A_2^*} + o(\eta^2),
\end{eqnarray}
where $A_{\eta}^*$ and $A_{per}^*$ are the homogenized tensors associated with $A_{\eta}$ and $A_{per}$ respectively, and the first and second-order corrections $\bar{A_1^*}$ and $\bar{A_2^*}$ can be, loosely speaking, computed in terms of
the microscopic properties of $A_{per}$ and $C_{per}$ and the statistics of second order of the random field $b_{\eta}$. The formulation is made precise in~\cite{ALB1} and in Sections~\ref{model}
and~\ref{heuristic} below.
\\

Motivations behind setting (\ref{problemintro}), as well as a review of the mathematical literature on similar issues and a comprehensive bibliography, can be found in \cite{ALB1}. We complement our study of the perturbative approach introduced with \cite{ALB1} in two different directions.\\

In Section \ref{model}, we rigorously establish an asymptotic
expansion of the homogenized tensor in a mathematical setting where 
our input parameter (the field $b_{\eta}$ in (\ref{problemintro})) enjoys
appropriate weak convergence properties, as $\eta$ vanishes, in a
reflexive Banach space, namely a Lebesgue space $L^\infty({\mathcal
  D},L^p(\Omega))$ (with $p>1$). In such a setting, we are in position
to rigorously prove a first order asymptotic expansion (announced in~\cite{ALBcras} and
precisely stated in~\cite[Théorème 2.1]{ALBcras} and Theorem~\ref{rig} below) for
the homogenization of $A_{\eta}$, using simple functional analysis
techniques very similar to those exposed in \cite{BLL}. In our
Corollaries~\ref{ordre2} and~\ref{correl}, the expansion is
pushed to second order under additional assumptions.\\

Our aim in Section \ref{heuristic} is to further extend our formal theory
of~\cite{ALB1}. Recall that this formal theory, rather than manipulating
the random field $b_{\eta}$ itself, consists in focusing on its law. We
indeed assume that the image measure (the law) corresponding
to the perturbation admits an expansion (see (\ref{pushexpand}) below) with respect to $\eta$ in the sense of  distributions.
 While \cite{ALB1}
has only addressed the specific case of a Bernoulli law, we consider here
more general laws and proceed with the same formal derivations. These derivations lead to a first-order correction $\bar{A_1^*}$ in (\ref{expansintro}) obtained as the limit when $N \rightarrow \infty$ of a sequence of tensors $A_{1}^{*,N}$ computed on the supercell $[-\frac{N}{2},\frac{N}{2}]^d$. It is the purpose of Proposition \ref{conva1} to prove the convergence of $A_{1}^{*,N}$. The second-order term $\bar{A_2^*}$ is likewise defined as a limit, up to extraction, of a sequence of tensors $A_{2}^{*,N}$ when $N \rightarrow \infty$. The proof of the boundedness of the sequence $A_{2}^{*,N}$ and thus of its convergence up to extraction is not given here for it involves long and technical computations. We refer to \cite{these_Arnaud} for the details. As in \cite{ALB1}, our approach in this Section exhibits close ties with classical defect-type theories used in solid state physics. \\ 

We emphasize that, in sharp contrast to the exact stochastic homogenization of $A_{\eta}$, the determination of the first and second-order terms in (\ref{expansintro}) 
relies on entirely deterministic computations, albeit of very different kind, for both approaches of Sections \ref{model} and \ref{heuristic}.\\   

Finally, a comprehensive series of numerical tests in Section \ref{num} show, beyond those contained in $\cite{ALB1}$, that the two approaches exposed here are efficient and quite robust: the computational workload induced by the perturbative approach is light compared to the direct homogenization of $\cite{ALB1}$, and expansion (\ref{expansintro}) proves to be accurate for not so small~perturbations.\\

We complement the text by a long appendix. The reader less interested in
theoretical issues can easily omit the reading of this appendix. Besides
providing, in Sections~\ref{distheo} and \ref{technical} and for consistency, some theoretical results useful
in the body of the text, the purpose of this appendix is two-fold. We examine in details
in Section~\ref{onedimension} 
the one-dimensional setting, and we show that, expectedly,  all our formal expansions
can be made rigorous through explicit computations. We next demonstrate, in Section~\ref{consistency},
that our two modes of derivation coincide
in a particular setting  appropriate for both the theoretical results of
Section~\ref{model} and the formal results of Section~\ref{heuristic}. This final
section therefore provides a \emph{proof} of our formal manipulations of
Section~\ref{heuristic}, in a setting -- we concede it -- that is not the setting
the  approach was designed to specifically address. Definite
conclusions on the theoretical validity of the approach developped in
Section~\ref{heuristic} are yet to be obtained, even though applicability and
efficiency are beyond doubt.\\

Throughout this paper, and unless otherwise mentioned, $C$ denotes a
constant that depends at most on the ambient dimension~$d$, and on the
tensors $A_{per}$ and $C_{per}$.  We write $C(\gamma)$ when $C$ depends
on $\gamma$ and possibly on $d$, $A_{per}$ and $C_{per}$. The
indices~$i$ and $j$ denote indices in $\llbracket 1, d
\rrbracket$.


\section{A model of a weakly randomly perturbed material} \label{model}

For consistency, we first recall the general setting of our related work \cite{ALB1}.\\

Throughout this article $\left(\Omega, \mathcal{F}, \mathbb{P} \right)$ denotes a probability space with $\mathbb{P}$ the probability measure and $\omega \in \Omega$ an event.  We denote by $\mathbb{E}(X)$ the expectation of a random variable $X$ and $Var(X)$ its variance.\\

We assume that the group $(\ZZ^d,+)$ acts on $\Omega$ and denote by $\tau_k, k \in \mathbb{Z}^d,$ the group action. We also assume that this action is measure-preserving, that is,
$$\forall \mathcal{A} \in \mathcal{F},  \forall k \in \mathbb{Z}^d, \; \mathbb{P}(\mathcal{A}) = \mathbb{P}(\tau_k \mathcal{A}),$$
and ergodic: 
$$\forall \mathcal{A} \in \mathcal{F},  (\forall k \in \mathbb{Z}^d, \mathcal{A} = \tau_k \mathcal{A})  \implies (\mathbb{P}(\mathcal{A})=0 \; \mathrm{or} \; \mathbb{P}(\mathcal{A})=1).$$

We call $F \in L^1_{loc}(\RR^d, L^1(\Omega))$ stationary if 
\begin{eqnarray}\label{statio}
 \forall k \in \ZZ^d, \; F(x+k, \omega) = F(x, \tau_k \omega)\; \; \mathrm{almost \; everywhere \; in \;} x \in \RR^d \; \mathrm{and \; } \omega \in \Omega.
\end{eqnarray}

Notice that if $F$ is deterministic, the notion of stationarity used here reduces to \\ $\ZZ^d$-periodicity, that is,
\begin{eqnarray}\label{periodi}
\forall k \in \ZZ^d, F(x+k) = F(x) \; \; \mathrm{almost \; everywhere \; in \;} x \in \RR^d.
\end{eqnarray}

We then consider the tensor field from $\RR^d \times \Omega$ to $\RR^{d \times d}$:
\begin{eqnarray} \label{homogeta}
A_{\eta}(x,\omega) = A_{per}(x) + b_{\eta}(x,\omega)C_{per}(x),
\end{eqnarray}

where $A_{per}$ and $C_{per}$ are two deterministic $\mathbb{Z}^d$-periodic tensor fields and $b_{\eta}$ a stationary ergodic scalar field. The matrix $A_{per}$ models the reference periodic material, perturbed by $C_{per}$. This perturbation is random, thus the presence of $b_{\eta}$.  We refer the reader to \cite{BLL} for a more detailed presentation of the stationary ergodic setting in a similar weakly random framework.\\

We make the following assumptions on the random field $b_{\eta}$: 
\begin{eqnarray} 
&\exists M >0, \forall \eta >0, \|b_{\eta}\|_{L^{\infty}(Q \times \Omega)} \leq M, \label{hyp1} \\
& \|b_{\eta}\|_{L^{\infty}(Q; L^{2}(\Omega))} \underset{\eta \rightarrow 0^{+}}{\rightarrow 0}, \label{hyp2}
\end{eqnarray}

where $Q$ is the unit cell $[-\frac{1}{2}, \frac{1}{2}]^d$.\\

Assumption (\ref{hyp2}) encodes that the perturbation for small $\eta$ is a rare event. Still, it is able to significantly modify the local structure of the material when it happens, for we do not require it to be small in $L^{\infty}(Q \times \Omega)$ as $\eta \rightarrow 0$.\\ 

We additionally assume that there exist $0 < \alpha \leq \beta$ such that for all $\xi \in \RR^d$, for almost all $x \in \RR^d$ and for all $s \in [-M,M]$,
\begin{eqnarray}\label{alphabound}
\alpha |\xi|^2 \leq A_{per}(x)\xi\cdot \xi, \qquad \alpha |\xi|^2 \leq \left(A_{per}+ s C_{per}\right)(x)\xi\cdot \xi,
\end{eqnarray}
\begin{eqnarray}\label{betabound}
A_{per}(x) \xi| \leq \beta |\xi|, \qquad |\left(A_{per}+ sC_{per}\right)(x)\xi| \leq \beta |\xi|.
\end{eqnarray}
We can therefore use the classical stochastic homogenization results (see for instance \cite{JKO} for a comprehensive review or \cite{ALB1} for a concise presentation). The cell problems associated with (\ref{homogeta}) read
\begin{eqnarray} \label{cellsol}
\left  \{
\begin{aligned}
& -\mathrm{div}\left(A_{\eta}(\nabla w_{i}^{\eta} +e_i ) \right) = 0 \quad \mathrm{in} \, \, \mathbb{R}^d, \\
& \nabla w_{i}^{\eta} \; \; \mathrm{stationary}, \; \; \mathbb{E}\left(\int_Q \nabla w_{i}^{\eta}\right) = 0.
\end{aligned}
\right.
\end{eqnarray}

Problem (\ref{cellsol}) has a solution unique up to the addition of a random constant. The function $w_{i}^{\eta}$ is called the $i$-th corrector or cell solution. \\

The homogenized tensor $A_{\eta}^{*}$ is given by
\begin{eqnarray} \label{homogal}
\forall i \in \llbracket 1,d \rrbracket, \; \; A_{\eta}^* e_i= \mathbb{E}\left(\int_Q A_{\eta}(\nabla w_{i}^{\eta} + e_i) \right).
\end{eqnarray}

Throughout the rest of this paper we will denote by $w_i^0$ the $i$-th cell solution associated with $A_{per}$, defined up to an additive constant by
\begin{eqnarray} \label{cellper}
\left  \{
\begin{aligned}
& -\mathrm{div}\left(A_{per}(\nabla w_{i}^{0} +e_i ) \right) = 0 \quad \mathrm{in} \, \, Q, \\
& w_{i}^{0} \; \; \ZZ^d-\mathrm{periodic}.
\end{aligned}
\right.
\end{eqnarray}

The periodic homogenized tensor is then given by
\begin{eqnarray} \label{homper}
\forall i \in \llbracket 1,d \rrbracket, \; \; A_{per}^* e_i= \int_Q A_{per}(\nabla w_{i}^{0} + e_i).
\end{eqnarray}

Due to the specific form of $A_{\eta}$, the following zero-order result can be easily proved. The proof is actually the same as that in Lemma 1 of \cite{ALB1}, which relies on the fact that $\|b_{\eta}\|_{L^{\infty}(Q; L^{2}(\Omega))}$ converges to $0$ as $\eta$ tends to $0$.

\begin{lemma} \label{facile}
When $\eta \rightarrow 0$, $A_{\eta}^* \rightarrow A_{per}^{*}$.  
\end{lemma}

Our goal is to find an asymptotic expansion for $A_{\eta}$ with respect to $\eta$, and a first answer is given by the following theorem announced as Théorème 1 in \cite{ALBcras} :
\begin{theorem}[Théorème 1, \cite{ALBcras}]  \label{rig}
Assume that $b_{\eta}$ satisfies (\ref{hyp1}) and (\ref{hyp2}), and denote by $m_{\eta} = \|b_{\eta}\|_{L^{\infty}(Q; L^{2}(\Omega))}$. There exists a subsequence of $\eta$, still denoted $\eta$ for the sake of simplicity, such that $\frac{b_{\eta}}{m_{\eta}}$ converges weakly-* in ${L^{\infty}(Q; L^{2}(\Omega))}$ to a limit field denoted by $\bar{b}_0$ when $n \rightarrow 0$. Then 
\begin{itemize}
\item for all $i \in \llbracket 1, d \rrbracket$, the following expansion
\begin{eqnarray} \label{expancell}
\nabla w_i^{\eta} = \nabla w_i^0 + m_{\eta} \nabla v_i^0 + o(m_{\eta})
\end{eqnarray}
holds weakly in $L^2(Q;L^2(\Omega))$, where $w_i^0$ is the solution to the $i$-th periodic cell problem and $v_i^{0}$ is solution to 
\begin{eqnarray} \label{vo}
\left  \{
\begin{aligned}
& -\mathrm{div}(A_{per} \nabla v_i^{0}) = \mathrm{div}\left(\bar{b}_0C_{per}(\nabla w^{0} + e_i)\right) \quad \mathrm{in} \, \, \mathbb{R}^d,\\
& \nabla v_i^{0} \; \; \mathrm{stationary}, \; \mathbb{E}\left(\int_Q \nabla v_i^{0}\right) = 0.
\end{aligned}
\right.
\end{eqnarray}
\item $A_{\eta}^*$ can be expanded up to first order as
\begin{eqnarray} \label{expans}
A_{\eta}^{*} = A_{per}^{*}+ m_{\eta} \tilde{A_1^{*}} + o(m_{\eta}),
\end{eqnarray}
where
\begin{eqnarray} \label{tilda1}
\forall i \in \llbracket 1, d \rrbracket, \; \tilde{A_1^{*}} e_i  = \int_Q  \mathbb{E}(\bar{b}_0)C_{per} (\nabla w_i^{0} +  e_i) + \int_Q A_{per}\nabla  \mathbb{E}(v_i^{0}).
\end{eqnarray}
\end{itemize}
\end{theorem}

\begin{proof}
We fix $i \in \llbracket 1, d \rrbracket$ and define $v_{i}^{\eta}= \frac{w_i^{\eta}-w_i^0}{m_{\eta}}$. $v_{i}^{\eta}$ is solution to
\begin{eqnarray} \label{vd} 
\left  \{
\begin{aligned}
& -\mathrm{div}\left(A_{\eta}\nabla v_{i}^{\eta} \right) = \mathrm{div}\left(\frac{b_{\eta}}{m_{\eta}} C_{per} \left(\nabla w_i^0 +e_i \right) \right) \quad \mathrm{in} \, \, \mathbb{R}^d, \\
 &\nabla v_{i}^{\eta} \; \; \mathrm{stationary}, \; \mathbb{E}\left(\int_Q \nabla v_{i}^{\eta}\right) = 0.
\end{aligned}
\right.
\end{eqnarray} 

Using an argument similar to that used in the proof of Lemma 1 in \cite{ALB1}, we have
$$\forall \eta > 0, \, \, \|\nabla v_i^{\eta}\|_{L^2(Q \times \Omega)} \leq \frac{1}{\alpha}\|C_{per} \left(\nabla w_i^0 +e_i \right)\|_{L^2(Q)}.$$
where $\alpha$ is defined in (\ref{alphabound}).\\

The sequence $\nabla v_{i}^{\eta}$ is bounded in $L^2(Q \times \Omega)$ and therefore, up to extraction, weakly converges in $L^2(Q \times \Omega)$ to some limit which is necessarily a gradient and which we denote $\nabla v_i^0$. Since $b_{\eta}$ converges strongly to $0$ in $L^2(Q \times \Omega)$, $b_{\eta} \nabla v_i^{\eta}$ converges to $0$ in $\mathcal{D}'(Q \times \Omega)$. It is then easy to pass to the limit $\eta \rightarrow 0$ in (\ref{vd}) and to deduce that $v_i^0$ is solution to 
\begin{eqnarray*} 
\left  \{
\begin{aligned}
& -\mathrm{div}\left(A_{per}\nabla v_{i}^{0} \right) = \mathrm{div}\left(\bar{b}_0 C_{per} \left(\nabla w_i^0 +e_i \right) \right) \quad \mathrm{in} \, \,\mathbb{R}^d, \\
 &\nabla v_{i}^{0} \; \; \mathrm{stationary}, \; \mathbb{E}\left(\int_Q \nabla v_{i}^{0}\right) = 0.
\end{aligned}
\right.
\end{eqnarray*}

Thus $\frac{\nabla w_i^{\eta}- \nabla w_i^0}{m_{\eta}}$ converges, up to extraction, weakly to $\nabla v_{i}^{0}$ in $L^2(Q \times \Omega)$. This amounts to say that we have the following first-order expansion:
$$ \nabla w_i^{\eta} = \nabla w_i^0 + m_{\eta} \nabla v_i^0 + o(m_{\eta}) \quad \mathrm{in} \; L^2(Q \times \Omega) \; \mathrm{weak}. $$

Inserting this expansion in (\ref{homogal}), we obtain
\begin{eqnarray*}
A_{\eta}^{*} e_i  &=& A_{per}^{*}e_i + m_{\eta}\int_Q  \mathbb{E}(\bar{b}_0)C_{per} (\nabla w_{i}^{0} +  e_i) 
+ m_{\eta}\int_Q A_{per}\nabla \mathbb{E}( v_{i}^{0}) + o(m_{\eta}), 
\end{eqnarray*}

which concludes the proof.

\end{proof}
\begin{remark}
Notice that taking the expectation of both sides of (\ref{vo}), $\mathbb{E}(v_{i}^{0})$ is actually the $\ZZ^d$-periodic function that is the unique solution (up to an additive constant) to 
\begin{eqnarray} \label{periodv}
\left  \{
\begin{aligned}
& -\mathrm{div}\left(A_{per}\nabla \mathbb{E}(v_{i}^{0}) \right) = \mathrm{div}\left(\mathbb{E}(\bar{b}_0) C_{per} \left(\nabla w_i^0 +e_i \right) \right) \quad \mathrm{in} \, \, Q,\\
 & \mathbb{E}(v_{i}^{0}) \; \; \mathbb{Z}^d-\mathrm{periodic}. \\
\end{aligned}
\right.
\end{eqnarray}
The computation of $A_{\eta}^*$ up to the first order in $m_{\eta}$ only requires solving $2d$ deterministic problems, namely (\ref{cellper}) and (\ref{periodv}), in the unit cell $Q$. 
\end{remark}

In fact, the situation is even more advantageous when $A_{per}$ is a symmetric matrix, as shown by our next remark.

\begin{remark} \label{trick}
Defining the adjoint problems to the cell problems (\ref{cellper}), 
\begin{eqnarray} \label{celladj}
\left  \{
\begin{aligned}
& -\mathrm{div}\left(A_{per}^T(\nabla \tilde{w}_{i}^0 + e_i ) \right) = 0 \quad \mathrm{in} \, \, Q,\\
& \tilde{w}_{i}^0 \; \; \mathbb{Z}^d-\mathrm{periodic}, 
\end{aligned}
\right.
\end{eqnarray}
where we have denoted by $A_{per}^T$ the transposed matrix of $A_{per}$, allows to write the first-order correction (\ref{tilda1}) in a slightly different form. Indeed, multiplying (\ref{periodv}) by $\tilde{w}_{j}^{0}$ and integrating by parts, we obtain 
$$ \int_Q A_{per}\nabla \mathbb{E}(v_{i}^{0}) \cdot \nabla \tilde{w}_{j}^{0}  = - \int_Q \mathbb{E}(\bar{b}_0) C_{per} \left(\nabla w_i^0 +e_i \right) \cdot \nabla \tilde{w}_{j}^{0}.$$
Likewise, multiplying (\ref{celladj}) by $\nabla \mathbb{E}(v_{i}^{0})$ and integrating by parts yields
$$ \int_Q A_{per} \nabla \mathbb{E}(v_{i}^{0}) \cdot \left(\nabla \tilde{w}_{j}^{0} +e_j \right) =0.$$
Combining these equalities gives
$$ \int_Q A_{per} \nabla \mathbb{E}(v_{i}^{0}) \cdot e_j = \int_Q \mathbb{E}(\bar{b}_0) C_{per} \left(\nabla w_i^0 +e_i \right) \cdot \nabla \tilde{w}_{j}^{0},$$
and thus (\ref{tilda1}) may be equivalently phrased as
\begin{eqnarray} \label{other}
\forall (i,j) \in \llbracket 1, d \rrbracket^2, \; \tilde{A_{1}^{*}} e_i \cdot e_j  &=& \int_Q  \mathbb{E}(\bar{b}_0)C_{per} (\nabla w_{i}^{0} +  e_i) \cdot (\nabla \tilde{w}_{j}^{0} +e_j). 
\end{eqnarray}
When $A_{per}$ is symmetric, $\tilde{w}_{j}^{0} = w_j^0$, and solving the periodic cell problems (\ref{cellper}) suffices to determine $A_{\eta}^{*}$ up to the first order in $m_{\eta}$.
\end{remark}

Pushing expansion (\ref{expans}) to second order requires more information on $b_{\eta}$:

\begin{cor} \label{ordre2}
Assume in addition to (\ref{hyp1}) and (\ref{hyp2}) that 
\begin{eqnarray}\label{hypordre2}
b_{\eta} = \eta \bar{b}_0 + \eta^2 \bar{r}_0 + o(\eta^2) \quad \mathrm{weakly}-^* \; \mathrm{in\; } L^{\infty}(Q; L^{2}(\Omega)).
\end{eqnarray}
Then 
\begin{itemize}
\item for all $i \in \llbracket 1, d \rrbracket$, the following expansion
\begin{eqnarray} \label{expanscell}
\nabla w_i^{\eta} = \nabla w_i^0 + \eta \nabla v_i^0 + \eta^2 \nabla z_i^0 + o(\eta^2)
\end{eqnarray}
holds weakly in $L^2(Q;L^2(\Omega))$, 
where $z_i^{0}$ is solution to 
\begin{eqnarray} \label{zo}
\left  \{
\begin{aligned}
& -\mathrm{div}(A_{per} \nabla z_i^{0}) = \mathrm{div}\left(\bar{r}_0C_{per}(\nabla w_i^{0} + e_i)\right) + \mathrm{div}\left(\bar{b}_0C_{per}\nabla v_i^0\right) \quad \mathrm{in} \, \, \mathbb{R}^d,\\
& \nabla z_i^{0} \; \; \mathrm{stationary}, \; \mathbb{E}\left(\int_Q \nabla z_i^{0}\right) = 0.
\end{aligned}
\right.
\end{eqnarray}
\item $A_{\eta}^*$ can be expanded up to second order as
\begin{eqnarray} \label{expans2}
A_{\eta}^{*} = A_{per}^{*} + \eta \tilde{A_1^*} + \eta^2 \tilde{A_2^*} + o(\eta^2), 
\end{eqnarray}
where $\tilde{A_1^*}$ is defined by (\ref{tilda1}) and for all $i \in \llbracket 1, d\rrbracket$,
\begin{eqnarray} \label{tilda2}
\tilde{A_2^*}e_i  =  \int_Q  \mathbb{E}(\bar{r}_0)C_{per} (\nabla w_i^{0} +  e_i) + \eta^2 \int_Q C_{per} \mathbb{E}(\bar{b}_0\nabla v_i^{0}) +\int_Q   A_{per} \nabla \mathbb{E}(z_i^{0}), 
\end{eqnarray}
\end{itemize}
or equivalently, for all $(i,j) \in \llbracket 1, d\rrbracket^2$,
\begin{eqnarray} \label{tilda2bis}
\tilde{A_2^*} e_i \cdot e_j  = \int_Q  \mathbb{E}(\bar{r}_0)C_{per} (\nabla w_i^{0} +  e_i) \cdot (\nabla \tilde{w_j}^0 + e_j) + \int_Q C_{per} \mathbb{E}(\bar{b}_0\nabla v_i^{0})\cdot (\nabla \tilde{w_j}^0 + e_j). 
\end{eqnarray}

\end{cor}

\begin{proof}
The proof follows the same pattern as that of Theorem \ref{rig}. The computation of the second order relies on the fact that (\ref{hypordre2}) implies that $\frac{b_{\eta}}{\eta}$ converges strongly to $\bar{b}_0$ in $L^{\infty}(Q; L^{2}(\Omega))$, whereas the convergence was weak in Theorem \ref{rig}. Likewise, the expansion of the cell solution, namely (\ref{expanscell}), implies that $\frac{\nabla w_i^{\eta} - \nabla w_i^{0}}{\eta}$ converges strongly to $\nabla v_i^{0}$ in $L^2(Q;L^2(\Omega))$. We then obtain (\ref{expans2}) and (\ref{tilda2}) by inserting (\ref{expanscell}) in (\ref{homogal}), and deduce (\ref{tilda2bis}) from (\ref{tilda2}) as in Remark \ref{trick}.
\end{proof}

The computation of $A_{\eta}^{*}$ up to the order $\eta^2$ is much more intricate than that up to the order $\eta$, for it requires determining $\mathbb{E}(\bar{b}_0\nabla v_i^{0})$. Computing the periodic deterministic function $\mathbb{E}(v_i^{0})$ solution to the simpler problem (\ref{periodv}) is not sufficient in general. We have to determine the stationary random field $v_i^{0}$ solution to (\ref{vo}) in $\RR^d$.\\

It turns out that in a particular, practically relevant setting, we may still avoid solving the random problem (\ref{vo}). This setting presents the additional advantage to provide insight on the influence of spatial correlation.

\begin{cor} \label{correl}
Assume that $b_{\eta}$ is uniform in each cell of $\mathbb{Z}^d$,  and writes
\begin{eqnarray}\label{specific}
b_{\eta}(x,\omega) = \sum_{k \in \mathbb{Z}^d} \mathds{1}_{Q+k}(x) B _{\eta}(\tau_k \omega),
\end{eqnarray}
where $B_{\eta}$ satisfies
\begin{eqnarray} 
&\forall \eta >0, \|B_{\eta}\|_{L^{\infty}(\Omega)} \leq M, \label{hyp1bis} \\
&B_{\eta} = \eta \bar{B}_0 + \eta^2 \bar{R}_0 + o(\eta^2) \quad \mathrm{weakly \; in\; } L^2(\Omega). \label{hyp2bis}
\end{eqnarray} 
Assume also that
\begin{eqnarray} \label{cov}
\sum_{k \in \ZZ^d} |cov(\bar{B}_0, \bar{B}_0 (\tau_k \cdot))| < \infty.
\end{eqnarray}

Then the second-order term (\ref{tilda2bis}) can be rewritten
\begin{eqnarray} \label{expans2correl}
\begin{aligned}
\tilde{A_2^*} e_i \cdot e_j =& \mathbb{E}(\bar{R}_0) \int_Q C_{per} (\nabla w_i^{0} +  e_i) \cdot (\nabla \tilde{w_j}^0 + e_j) + 
 Var(\bar{B}_0) \int_{Q} C_{per} \nabla t_i \cdot (\nabla \tilde{w}_j^0 + e_j)  \\ 
&+ (\mathbb{E}(\bar{B}_0))^2 \int_Q C_{per} \nabla s_i  \cdot (\nabla \tilde{w}_j^0 + e_j) \\
& + \sum_{k \in \ZZ^d, k\neq 0} cov(\bar{B}_0, \bar{B}_0 (\tau_k \cdot)) \int_Q C_{per} \nabla t_i( \cdot -k) \cdot (\nabla \tilde{w}_j^0 + e_j),
\end{aligned}
\end{eqnarray}
where $t_i$ is a $L^2_{loc}(\RR^d)$ function solving
\begin{eqnarray} \label{ti} 
\left  \{
\begin{aligned}
& -\mathrm{div}\left(A_{per}\nabla t_{i} \right) = \mathrm{div}\left(C_{per} \mathds{1}_{Q} \left(\nabla w_i^0 +e_i \right) \right) \quad \mathrm{in} \, \, \mathbb{R}^d, \\
 &\nabla t_{i} \in L^2(\RR^d),
\end{aligned}
\right.
\end{eqnarray}
and $s_i$ solves
\begin{eqnarray} \label{si}
\left  \{
\begin{aligned}
& -\mathrm{div}\left(A_{per}\nabla s_{i} \right) = \mathrm{div}\left(C_{per}\left(\nabla w_i^0 +e_i \right) \right) \quad \mathrm{in} \, \, Q,\\
 & s_i \; \; \ZZ^d-\mathrm{periodic}. 
\end{aligned}
\right.
\end{eqnarray}

\end{cor}

\begin{proof}

We notice first that the specific form (\ref{specific}) of $b_{\eta}$ considered implies that $\bar{b}_0$ and $\bar{r}_0$ defined in (\ref{hypordre2}) here write
\begin{eqnarray}\label{specificb0}
\bar{b}_0(x,\omega) = \sum_{k \in \mathbb{Z}^d} \mathds{1}_{Q+k}(x) \bar{B}_0(\tau_k\omega),
\end{eqnarray}
\begin{eqnarray}\label{specificr0}
\bar{r}_0(x,\omega)  = \sum_{k \in \mathbb{Z}^d} \mathds{1}_{Q+k}(x) \bar{R}_0(\tau_k\omega).
\end{eqnarray}

The rest of the proof mainly consists in showing that in this particular setting, $\nabla v_i^0$ and the product $\bar{b}_0 \nabla v_i^0$ can be written using the deterministic functions $t_i$ and $s_i$. The existence of $t_i$ and its uniqueness up to an additive constant come from Lemmas 6 and 7 in \cite{ALB1}.\\

  We start by proving that the sum 
\begin{eqnarray} \label{series}
\sum_{k \in \mathbb{Z}^d}\left(\bar{B}_0(\tau_k\omega) - \mathbb{E}(\bar{B}_0)\right) \nabla t_i(x-k)
\end{eqnarray}
is a convergent series in $L^2(Q \times \Omega)$.\\

To this end, we compute the norm of the remainder of this series:
\begin{eqnarray*}
\begin{aligned}
& \left\|\sum_{|k| \geq N}\left(\bar{B}_0(\tau_k \cdot) - \mathbb{E}(\bar{B}_0)\right) \nabla t_i(\cdot - k) \right\|^2_{L^2(Q \times \Omega)} \\
& = \sum_{|k| \geq N} \sum_{|l| \geq N} cov(\bar{B}_0(\tau_k \cdot), \bar{B}_0 (\tau_l \cdot)) \int_Q \nabla t_i(\cdot - k) \nabla t_i(\cdot - l) \\
& \leq \frac{1}{2} \sum_{|k| \geq N} \sum_{|l| \geq N} |cov(\bar{B}_0(\tau_k \cdot), \bar{B}_0 (\tau_l \cdot))| (\|\nabla t_i(\cdot - k)\|_{L^2(Q)}^2 + \|\nabla t_i(\cdot - l)\|_{L^2(Q)}^2) \\
& \leq \sum_{|k| \geq N} \sum_{|l| \geq N} |cov(\bar{B}_0(\tau_k \cdot), \bar{B}_0 (\tau_l \cdot))| \, \|\nabla t_i(\cdot - k)\|_{L^2(Q)}^2 \\
& \leq \sum_{|k| \geq N} \left( \|\nabla t_i(\cdot - k)\|_{L^2(Q)}^2  \sum_{|l| \geq N} |cov(\bar{B}_0(\tau_k \cdot), \bar{B}_0 (\tau_l \cdot))| \right) \\
& \leq \sum_{|k| \geq N} \left( \|\nabla t_i(\cdot - k)\|_{L^2(Q)}^2  \sum_{|l| \geq N} |cov(\bar{B}_0, \bar{B}_0 (\tau_{l-k} \cdot))| \right) \\
& \leq  \sum_{|k| \geq N} \|\nabla t_i(\cdot - k)\|_{L^2(Q)}^2  \sum_{k \in \ZZ^d} |cov(\bar{B}_0, \bar{B}_0 (\tau_{k} \cdot))|.
\end{aligned}
\end{eqnarray*}

Using (\ref{cov}), we obtain
\begin{eqnarray} \label{borncov}
\left\|\sum_{|k| \geq N}\left(\bar{B}_0(\tau_k \cdot) - \mathbb{E}(\bar{B}_0)\right) \nabla t_i(\cdot - k) \right\|^2_{L^2(Q \times \Omega)} \leq C  \sum_{|k| \geq N} \|\nabla t_i(\cdot - k)\|_{L^2(Q)}^2.
\end{eqnarray}

Since $\nabla t_i \in L^2(\RR^d)$, the right-hand side of (\ref{borncov}) converges to zero when $N$ goes to infinity.\\

Consequently, (\ref{series}) defines a vector $T$ in $L^2(Q \times \Omega)$. It is clear from (\ref{series}) that $\frac{\partial T_p}{\partial x_n} = \frac{\partial T_n}{\partial x_p}$ for all $(n,p) \in \llbracket 1, d \rrbracket^2$. Thus $T$ is a gradient, and there exists a function $\tilde{v}_i$ such that
\begin{eqnarray} \label{v_i}
\nabla \tilde{v}_i = T + \mathbb{E}(\bar{B}_0) \nabla s_i = \sum_{k \in \mathbb{Z}^d}\left(\bar{B}_0(\tau_k \cdot) - \mathbb{E}(\bar{B}_0)\right) \nabla t_i(x-k) + \mathbb{E}(\bar{B}_0) \nabla s_i. 
\end{eqnarray}

Since $s_i$ is $\ZZ^d$-periodic, we deduce from (\ref{v_i}) that
\begin{eqnarray} \label{bconditions}
\nabla \tilde{v}_i \; \mathrm{is \; stationary \; and \; } \mathbb{E}\left(\int_Q \nabla \tilde{v}_i \right) = 0.    
\end{eqnarray}

We then compute, using (\ref{ti}) and (\ref{si}),
\begin{eqnarray} \label{tedious}
-\mathrm{div} (A_{per} \nabla  \tilde{v}_i) &= &\sum_{k \in \mathbb{Z}^d} -\mathrm{div} (A_{per} \nabla  t_i(\cdot -k)) \left(\bar{B}_0(\tau_k \cdot) - \mathbb{E}(\bar{B}_0)\right) \nonumber \\
&& -\mathrm{div} (A_{per} \nabla  s_i)\mathbb{E}(\bar{B}_0) \nonumber \\
&=& \sum_{k \in \mathbb{Z}^d} \mathrm{div}\left(C_{per} \mathds{1}_{Q+k} \left(\nabla w_i^0 +e_i \right) \right) \left(\bar{B}_0(\tau_k \cdot) - \mathbb{E}(\bar{B}_0)\right) \nonumber \\
&& + \mathrm{div}\left(C_{per}\left(\nabla w_i^0 +e_i \right) \right) \mathbb{E}(\bar{B}_0) \nonumber \\
&=& \sum_{k \in \mathbb{Z}^d} \mathrm{div}\left(C_{per} \mathds{1}_{Q+k} \, \bar{B}_0(\tau_k \cdot) \left(\nabla w_i^0 +e_i \right) \right). 
\end{eqnarray}

Because of (\ref{specificb0}), (\ref{tedious}) implies 
\begin{eqnarray} \label{equality}
 -\mathrm{div} (A_{per} \nabla  \tilde{v}_i) =  \mathrm{div}\left(\bar{b}_0C_{per}(\nabla w^{0} + e_i)\right).
\end{eqnarray}

It follows from (\ref{bconditions}) and (\ref{equality}) that $\tilde{v}_i$ solves (\ref{vo}). As (\ref{vo}) has a solution unique up to the addition of a random constant, we obtain
\begin{eqnarray}\label{specificv0}
\nabla v_i^0 = \nabla \tilde{v}_i = \sum_{k \in \mathbb{Z}^d}\left(\bar{B}_0(\tau_k\cdot) - \mathbb{E}(\bar{B}_0)\right) \nabla t_i(x-k) + \mathbb{E}(\bar{B}_0) \nabla s_i.
\end{eqnarray}

We deduce from (\ref{specificb0}) and (\ref{specificv0}) that
\begin{eqnarray*}
\mathbb{E}(\bar{b}_0 \nabla v_i^0) &=& \sum_{k \in \ZZ^d} \sum_{l \in \ZZ^d } \mathds{1}_{Q+l} \, \mathbb{E}(\bar{B}_0(\tau_l \cdot) (\bar{B}_0(\tau_k \cdot) - \mathbb{E}(\bar{B}_0)))\nabla t_i(\cdot-k)\\
 & &+ (\mathbb{E}(\bar{B}_0))^2 \sum_{l \in \ZZ^d } \mathds{1}_{Q+l} \nabla s_i \\
&=& \sum_{k \in \ZZ^d} \sum_{l \in \ZZ^d } \mathds{1}_{Q+l} \, cov(\bar{B}_0(\tau_k \cdot), \bar{B}_0 (\tau_l \cdot)) \nabla t_i(\cdot-k)  + (\mathbb{E}(\bar{B}_0))^2 \sum_{l \in \ZZ^d } \mathds{1}_{Q+l} \nabla s_i,
\end{eqnarray*}
and then that
\begin{eqnarray} \label{specificbovo}
\begin{aligned}
\mathds{1}_Q \mathbb{E}(\bar{b}_0 \nabla v_i^0) =& Var(\bar{B}_0)\nabla t_i + \sum_{k \in \ZZ^d, k \neq 0} cov(\bar{B}_0(\cdot), \bar{B}_0 (\tau_k \cdot))\nabla t_i(\cdot-k) \\ 
& + (\mathbb{E}(\bar{B}_0))^2 \nabla s_i.
\end{aligned}
\end{eqnarray}

We conclude by inserting (\ref{specificr0}) and (\ref{specificbovo}) in (\ref{tilda2bis}).

\end{proof}

Theorem \ref{rig} (and its two corollaries) are only of interest if $\mathbb{E}(\bar{b}_0) \neq 0$. Indeed, if $\mathbb{E}(\bar{b}_0)=0$ it only states that $A_{\eta}^* = A_{per}^* + o(m_{\eta})$.\\

The prototypical case where Theorem \ref{rig} does not provide valuable information is the case studied in \cite{ALB1}: $b_{\eta}(x,\omega) = \sum_{k \in \mathbb{Z}^d} \mathds{1}_{Q+k}(x) B _{\eta}^k(\omega)$, where the $B _{\eta}^k$ are independent identically distributed variables that have Bernoulli law with parameter $\eta$, i.e are equal to 1 with probability $\eta$ and to 0 with probability $1-\eta$. Then, using the notation of Theorem~\ref{rig}, $b_{\eta}^2=b_{\eta}$, $m_{\eta} = \sqrt{\eta}$ and  $\bar{b}_0=0$, and we only get $A_{\eta}^* = A_{per}^* + o(\sqrt{\eta})$ (while appendix 6.1 of \cite{ALB1} shows that there exists a tensor $\bar{A_1^*}$ such that $A_{\eta}^* = A_{per}^* + \eta \bar{A_1^*} + o(\eta)$ at least in dimension one). Omitting the dependence on the space variables since $b_{\eta}$ is uniform in each cell of $\ZZ^d$ in this particular setting, a suitable functional space $F$ on $\Omega$ to obtain a non trivial weak limit of $\frac{b_{\eta}}{\|b_{\eta}\|_F}$ would be $L^1(\Omega)$ for the norm of each $B _{\eta}^k$ in $L^1(\Omega)$ is equal to $\eta$. The Dunford-Petti weak compactness criterion in that space is however not satisfied by $\frac{b_{\eta}}{\|b_{\eta}\|_{L^1(\Omega)}}$. The reason is of course that $\frac{b_{\eta}}{\|b_{\eta}\|_{L^1(\Omega)}}$ converges in the set of bounded measures to a Dirac mass. The techniques used in the proof of Theorem \ref{rig} and its two corollaries thus do not work in this setting.\\  

The above considerations somehow suggest that an alternative viewpoint might be useful. Because of (\ref{hyp2}), the image measure $dP_{\eta}^x$ of $b_{\eta}(x,\cdot)$ converges to a Dirac mass in the sense of distributions. Our alternate approach, related to our work \cite{ALB1}, consists in working out an expansion of the image measure (or of the law), rather than an expansion of the random variable. Like in \cite{ALB1}, our manipulations are mostly formal. Some rigorous foundations, in specific settings, are provided in the appendix.

\section{A formal approach} \label{heuristic}

\subsection{A new assumption on the image measure}

For simplicity, we assume as in Corollary \ref{correl} that $b_{\eta}$ is uniform in each cell of $\mathbb{Z}^d$, and is of the form
\begin{eqnarray} \label{forme}
b_{\eta}(x,\omega) = \sum_{k \in \mathbb{Z}^d} \mathds{1}_{Q+k}(x) B _{\eta}^k(\omega),
\end{eqnarray}
where the $B _{\eta}^k$ are independent identically distributed random variables, the distribution of which is given by a "mother variable" $B_{\eta}$. For convenience we slightly modify (\ref{hyp1bis}) and require 
\begin{eqnarray} 
& \exists \epsilon > 0, \forall \eta >0, \|B_{\eta}\|_{L^{\infty}(\Omega)} \leq M - \epsilon \label{hyp1bise}\\ 
& \|B_{\eta}\|_{L^2(\Omega)} \underset{\eta \rightarrow 0}{\rightarrow 0}. \label{hyp2bise}
\end{eqnarray}

Assumption (\ref{hyp1bise}) is a technical assumption which implies in particular that for every $\eta >0$, the image measure $dP_{\eta}$ of $B_{\eta}$ is a distribution with compact support contained in the \emph{open} set $]-M,M[$. Of course the specific values of $M$ and $\epsilon$ have no particular significance. Throughout the sequel we denote by $\mathcal{E}'(]-M,M[)$ the space of distributions on $\RR$ with compact support in $]-M,M[$, and by $\langle T, \varphi \rangle$ the action of a distribution $T \in \mathcal{E}'(]-M,M[)$ on a test function $\varphi \in \mathcal{C}^{\infty}(]-M,M[)$. Basic elements of distribution theory are recalled in Section \ref{distheo} of the appendix, for convenience of the reader not familiar with technical issues.\\

Because of assumption (\ref{hyp2bise}) and Lebesgue dominated convergence theorem, it is clear that for every $\varphi \in \mathcal{C}^{\infty}(]-M,M[)$,
$$ \mathbb{E}(\varphi(B_{\eta}))  \underset{\eta \rightarrow 0^{+}}{\rightarrow} \varphi(0).$$
Since $\mathbb{E}(\varphi(B_{\eta})) = \langle dP_{\eta}, \varphi \rangle \; \mathrm{and} \; \varphi(0) =\langle \delta_0, \varphi \rangle$
where $\delta_0$ is the Dirac mass at $0$, $dP_{\eta}$ converges to $\delta_0$ in $\mathcal{E}'(]-M,M[)$.\\

This leads us to assume that $dP_{\eta}$ satisfies 
\begin{eqnarray} \label{pushexpand}
dP_{\eta} = \delta_0 + \eta d\bar{P}_1 + \eta^2 d\bar{P}_2 + o(\eta^2) \quad \mathrm{in} \; \mathcal{E}'(]-M,M[),
\end{eqnarray}
which is equivalent to
$$\forall \varphi \in \mathcal{C}^{\infty}(]-M,M[), \; \mathbb{E}(\varphi(B_{\eta})) = \langle dP_{\eta}, \varphi \rangle = \varphi(0) +  \eta \langle d\bar{P}_1, \varphi \rangle + \eta^2 \langle d\bar{P}_2, \varphi \rangle + o(\eta^2).$$

Of course $d\bar{P}_1$ and $d\bar{P}_2$ also have a compact support contained in $]-M,M[$\,: for every test function $\varphi$ with compact support in $\RR \backslash [-M+\epsilon,M-\epsilon]$, it holds for all $\eta >0$
$$\langle dP_{\eta}, \varphi \rangle =\mathbb{E}(\varphi(B_{\eta})) = 0 = \eta \langle d\bar{P}_1, \varphi \rangle + \eta^2 \langle d\bar{P}_2, \varphi \rangle + o(\eta^2),$$
which yields $\langle d\bar{P}_1, \varphi \rangle = \langle d\bar{P}_2, \varphi \rangle = 0$. Then the supports of $d\bar{P}_1$ and $d\bar{P}_2$ are contained in $[-M+\epsilon,M-\epsilon] \subset ]-M,M[$.\\

Denoting by $M' = M -\epsilon/2$, we deduce from Proposition \ref{finiteorder} of the appendix that there exists a constant $C>0$ and integers $p_1$ and $p_2$ (namely the orders of $d\bar{P}_1$ and $d\bar{P}_2$ respectively) such that
\begin{eqnarray} \label{major1}
\forall \varphi \in C^{\infty}(]-M,M[), \; |\langle d\bar{P}_1, \varphi \rangle| \leq C \sup_{s \in [-M',M']} \sup_{0 \leq n \leq p_1} \left|\frac{d^n}{ds^n} \varphi(s) \right|,
\end{eqnarray}

\begin{eqnarray} \label{major2}
\forall \varphi \in C^{\infty}(]-M,M[), \; |\langle d\bar{P}_2, \varphi \rangle| \leq C \sup_{x \in [-M',M']} \sup_{0 \leq n \leq p_2} \left|\frac{d^n}{ds^n} \varphi(s) \right|.
\end{eqnarray}

Let us now give some additional motivations underlying assumption (\ref{pushexpand}).\\

The first motivation is related to our work presented in \cite{ALB1} in which $B_{\eta}$ has Bernoulli law with parameter $\eta$, meaning that it is equal to $1$ with probability $\eta$ and $0$ with probability $1-\eta$. Then the image measure $dP_{\eta}$ is equal to $\delta_0 + \eta (\delta_1 - \delta_0)$, so that it satisfies (\ref{pushexpand}) exactly at order $1$ with $d\bar{P}_1 = \delta_1 - \delta_0$.\\

The second motivation comes from the following result, which shows that there is an easy way, used in our numerical experiments, to build perturbations satisfying (\ref{pushexpand}).

\begin{lemma} \label{lemexp}
Consider $B$ a random variable in $L^{3}(\Omega)$. Let $K$ be a positive real, and define $B_{\eta} = \eta B \mathds{1}_{|\eta B| \leq K}$. Then $B_{\eta}$, which obviously satisfies (\ref{hyp1bise}) and (\ref{hyp2bise}), also satisfies (\ref{pushexpand}) with
\begin{eqnarray} \label{expanspart}
dP_{\eta} = \delta_0 - \eta \mathbb{E}(B)\delta_0'+ \frac{\eta^2}{2} \mathbb{E}(B^2) \delta_0'' + \mathcal{O}(\eta^3) \; \mathrm{in} \; \mathcal{E}'(\mathbb{R}).
\end{eqnarray}
\end{lemma}

\begin{proof}
Let us denote by $dP$ the image measure of $B$, and consider $\varphi \in \mathcal{D}(\RR)$ (i.e $\varphi \in \mathcal{C}^{\infty}(\RR)$ and has compact support). Then 
\begin{eqnarray}
\langle dP_{\eta}, \varphi \rangle  &=& \int_{|\eta s| \leq K} \varphi(\eta s) dP + \varphi(0) \int_{|\eta s| \geq K} dP \\
&=& \int_{\mathbb{R}} \varphi(\eta s) dP + \int_{|\eta s| \geq K} (\varphi(0) - \varphi(\eta s)) dP. 
\end{eqnarray}

Since $B$ is in $L^{3}(\Omega)$, 
$$ \int_{|\eta s| \geq K} dP = \mathcal{O}(\eta^{3}),$$
and thus, $\varphi$ being a bounded function,  $$\langle dP_{\eta}, \varphi \rangle  = \int_{\mathbb{R}} \varphi(\eta s) dP + \mathcal{O}(\eta^{3}).$$  Then, since $\varphi \in \mathcal{D}(\mathbb{R})$, there exists $C >0$ such that
$$ \forall s \in \mathbb{R}, \quad \left|\frac{\varphi(\eta s)-\varphi(0)-\eta s \varphi'(0) -\frac{\eta^2}{2} s^2\varphi''(0) }{\eta^2}\right| \leq C \eta |s|^3.$$

Again using $B \in L^{3}(\Omega)$, this implies that 
$$\int_{\mathbb{R}}  \left(\frac{\varphi(\eta s)-\varphi(0)-\eta s \varphi'(0) -\frac{\eta^2}{2} s^2\varphi''(0) }{\eta^2}\right) dP \underset{\eta \rightarrow 0}{\rightarrow 0}$$
which is just a rewriting of (\ref{expanspart}) since $\displaystyle \int dP =1$, $\displaystyle \int s dP = \mathbb{E}(B)$ and $\displaystyle \int s^2 dP = \mathbb{E}(B^2)$.
\end{proof}

Before exposing our approach in this new setting, we prove the following elementary result which we will often use in the sequel:

\begin{lemma} \label{distreasy}
It holds $\langle d\bar{P}_1, 1 \rangle = 0$ and $\langle d\bar{P}_2, 1 \rangle =0$. 
\end{lemma}

\begin{proof}
 It holds on the one hand $\langle dP_{\eta}, 1 \rangle = 1$ since $dP_{\eta}$ is a probability measure, and on the other hand
\begin{eqnarray*}
\langle dP_{\eta}, 1 \rangle &=& \langle \delta_0, 1 \rangle + \eta \langle d\bar{P}_1, 1 \rangle + \eta^2 \langle d\bar{P}_2, 1 \rangle + o(\eta^2) \\
&=&  1 + \eta \langle d\bar{P}_1, 1 \rangle + \eta^2 \langle d\bar{P}_2, 1 \rangle + o(\eta^2),
\end{eqnarray*}
so that the conclusion follows.

\end{proof}

\subsection{An ergodic approximation of the homogenized tensor}

Let us consider a specific realization $\tilde{\omega} \in \Omega$ of $A_{\eta}$ in $I_N = [-\frac{N}{2},\frac{N}{2}]^d$, $N$ being for simplicity an odd integer, and solve the following ``supercell'' problem:
\begin{eqnarray} \label{ergodiccorr} 
\left  \{
\begin{aligned}
& -\mathrm{div}\left(A_{\eta}(x,\tilde{\omega})(\nabla w_{i}^{\eta,N,\tilde{\omega}} + e_i ) \right) = 0 \quad \, \mathrm{in}  \, I_N,\\
& w_{i}^{\eta,N,\tilde{\omega}} \,(N\ZZ)^d-\mathrm{periodic}.
\end{aligned}
\right.
\end{eqnarray} 

Then we have 
\begin{eqnarray} \label{expect}
\forall i \in \llbracket 1, d \rrbracket, \; \; A_{\eta}^*e_i = \lim_{N \rightarrow + \infty} \frac{1}{N^d} \mathbb{E} \left(\int_{I_N}A_{\eta}(x,\omega)(\nabla w_{i}^{\eta,N,\omega}(x) + e_i)\right) dx. 
\end{eqnarray}

The proof of (\ref{expect}) is given in \cite{ALB1}. We only outline it here for convenience. We know from Theorem 1 in \cite{BP} that 
\begin{eqnarray} \label{almost}
\displaystyle \frac{1}{N^d} \int_{I_N}A_{\eta}(x,\tilde{\omega})(\nabla w_{i}^{\eta,N,\tilde{\omega}}(x) + e_i ) \, dx \mathrm{\; converges \; to\;} A_{\eta}^*e_i \mathrm{\; almost \; surely \; in} \; \tilde{\omega} \in \Omega.
\end{eqnarray} Since $\displaystyle \frac{1}{N^d} \int_{I_N}A_{\eta}(x,\tilde{\omega})(\nabla w_{i}^{\eta,N,\tilde{\omega}}(x) + e_i ) dx$ is the periodic homogenization of $A_{\eta}(x,\tilde{\omega})$ on $I_N$, it is also well known that for all $(i,j) \in \llbracket 1, d \rrbracket^2$,
\begin{eqnarray} 
\begin{aligned}
\frac{1}{N^d} \left(\int_{I_N}A_{\eta}^{-1}(x,\tilde{\omega})dx\right)^{-1} e_i \cdot e_j  \leq \frac{1}{N^d} \int_{I_N}A_{\eta}(x,\tilde{\omega})(\nabla w_{i}^{\eta,N,\tilde{\omega}}(x) + e_i )\cdot e_j dx \\
 \leq \frac{1}{N^d} \left(\int_{I_N}A_{\eta}(x,\tilde{\omega})dx\right) e_i \cdot e_j,
\end{aligned}
\end{eqnarray}

 so that for all $N \in 2\NN+1$, for all $\eta>0$ and for almost all $\tilde{\omega} \in \Omega$,
\begin{eqnarray} \label{Lebesgue}
\left|\frac{1}{N^d} \int_{I_N}A_{\eta}(x,\tilde{\omega})(\nabla w_{i}^{\eta,N,\tilde{\omega}}(x) + e_i )\cdot e_j dx \right| \leq \beta,
\end{eqnarray}
where $\beta$ is defined by (\ref{betabound}).
Using (\ref{Lebesgue}) and the Lebesgue dominated convergence theorem, we can take the expectation in (\ref{almost}) and get (\ref{expect}).

\begin{remark}
The same result holds for homogeneous Dirichlet and Neumann boundary conditions instead of periodic conditions in the definition of $w_{i}^{\eta,N,\tilde{\omega}}$ (see \cite{BP} for more details).
\end{remark}

For convenience, we label the unit cells of $I_N$ from $1$ to $N^d$. The $k$-th cell is denoted by $Q_k$, for $1\leq k \leq N^d$.
A given realization $A_{\eta}(x,\tilde{\omega})$ can then be rewritten
$$ A_{\eta}(x,\tilde{\omega}) = A_{per}(x) +  \sum_{k=1}^{N^d} \mathds{1}_{Q_k}(x) s_k C_{per}(x),$$
with $s_k = B_{\eta}^k(\tilde{\omega})$ for all $k \in \llbracket 1, N^d \rrbracket.$
The $B_{\eta}^k(\tilde{\omega})$ being independent random variables, the joint probability of the $N^d$-uplet $(s_1,\cdots,s_{N^d})$ is simply the product $\displaystyle  \prod_{k=1}^{N^d} dP_{\eta}(s_k)$.

\begin{remark}
The approach exposed in the sequel works also, with minor changes, for random variables which are not independent but correlated with a finite length of correlation. We present it in the independent setting for simplicity. 
\end{remark}

We now define $ \displaystyle A^{s_1, \cdots, s_{N^d}} = A_{per} +  \sum_{k=1}^{N^d} \mathds{1}_{Q_k} s_k C_{per}$ for $(s_1,\cdots,s_{N^d}) \in [-M,M]^{N^d}$. We denote by $w_i^{s_1, \cdots, s_{N^d}}$ the solution of the $i$-th cell problem for the periodic homogenization of $A^{s_1, \cdots, s_{N^d}}$ on $I_N$, that is
\begin{eqnarray}\label{sgen} 
\left  \{
\begin{aligned}
& -\mathrm{div}\left(A^{s_1, \cdots, s_{N^d}} (\nabla w_i^{s_1, \cdots, s_{N^d}} + e_i ) \right) = 0 \quad \, \mathrm{in}  \, I_N,\\
& w_i^{s_1, \cdots, s_{N^d}} \,(N\ZZ)^d-\mathrm{periodic}.
\end{aligned}
\right.
\end{eqnarray} 
Then,
defining
\begin{eqnarray} \label{aetastarn}
 A_{\eta}^{*,N} e_i= \frac{1}{N^d} \mathbb{E} \left(\int_{I_N}A_{\eta}(x,\omega)(\nabla w_{i}^{\eta,N,\omega}(x) + e_i)\right) dx,
\end{eqnarray}
we have 
\begin{eqnarray}\label{product}
 A_{\eta}^{*,N}e_i =  \frac{1}{N^d}\int_{\RR^{N^d}}\left(\int_{I_N} A^{s_1, \cdots, s_{N^d}} ( \nabla w_i^{s_1, \cdots, s_{N^d}} + e_i ) \right) \prod_{k=1}^{N^d} dP_{\eta}(s_k).
\end{eqnarray}
It is proved in Lemma \ref{regular} of the Appendix that $ \nabla w_i^{s_1, \cdots, s_{N^d}}$ is a $\mathcal{C}^{\infty}$ function of $(s_1, \cdots, s_{N^d})$ in $]-M,M[^{N^d}$. Thus, since $d\bar{P}_1$ and $d\bar{P}_2$ have compact support in $]-M,M[$ (as well as $\delta_0$ of course), we can make these distributions act on $A^{s_1, \cdots, s_{N^d}}$ and $\nabla w_i^{s_1, \cdots, s_{N^d}}$ as functions of $(s_1, \cdots, s_{N^d})$.\\

It follows from (\ref{pushexpand}) that 
\begin{eqnarray} \label{devprod}
\begin{aligned}
\prod_{k=1}^{N^d} dP_{\eta}(s_k) =& \prod_{k=1}^{N^d} \delta_0(s_k) + \eta \sum_{l=1}^{N^d} d\bar{P}_1(s_l)\prod_{k=1, k\neq l}^{N^d} \delta_0(s_k) \\
 &+ \frac{\eta^2}{2} \sum_{l=1}^{N^d} \sum_{m=1}^{N^d} d\bar{P}_1(s_l) d\bar{P}_1(s_m)\prod_{k=1, k\neq \{l,m\}}^{N^d} \delta_0(s_k)\\
 &+ \eta^2 \sum_{l=1}^{N^d} d\bar{P}_2(s_l)\prod_{k=1, k\neq l}^{N^d} \delta_0(s_k) + o_N(\eta^2) \; \; \mathrm{in} \, \mathcal{E}'(]-M,M[^{N^d}).
\end{aligned}
\end{eqnarray}

We stress that the remainder $o_N(\eta^2)$ in (\ref{devprod}) depends on $N$, hence the notation.\\

 Moreover the products (\ref{devprod}) are to be understood as tensorized products: we work in $\mathcal{E}'(]-M,M[) \otimes_1 \mathcal{E}'(]-M,M[) \otimes_2 \cdots   \otimes_{N^d-1} \mathcal{E}'(]-M,M[) \subset \mathcal{E}'(]-M,M[^{N^d})$.\\ 

Inserting (\ref{devprod}) in (\ref{product}), we obtain the following second-order expansion
\begin{eqnarray}\label{asymptN}
A_{\eta}^{*,N} = A_{0}^{*,N} + \eta A_{1}^{*,N} + \eta^2 A_{2}^{*,N} +o_N(\eta^2). 
\end{eqnarray}

Before making the first three orders in (\ref{asymptN}) precise, note that (\ref{expect}), (\ref{aetastarn}) and (\ref{asymptN}) imply
\begin{eqnarray} \label{realasympt}
A_{\eta}^{*} = \lim_{N \rightarrow \infty} \left(A_{0}^{*,N} + \eta A_{1}^{*,N} + \eta^2 A_{2}^{*,N} +o_N(\eta^2) \right)
\end{eqnarray}

In the sequel we exchange in (\ref{realasympt}) the limit in $N$ and the series in $\eta$  in order to guess a second-order expansion of $A_{\eta}^{*}$ depending only on $\eta$. Since we are not able to justify this permutation, our approach is formal.\\

We now detail the first three orders in (\ref{asymptN}).\\

First, we notice that for $i \in \llbracket 1, d \rrbracket$,
\begin{eqnarray*}
 A_{0}^{*,N} e_i &=& \frac{1}{N^d}  \left \langle \prod_{k=1}^{N^d} \delta_0(s_k),  \int_{I_N} A^{s_1, \cdots, s_{N^d}} ( \nabla w_i^{s_1, \cdots, s_{N^d}} + e_i )   \right \rangle \\
&=& \frac{1}{N^d}  \int_{I_N} A^{0, \cdots, 0} (\nabla w_i^{0, \cdots, 0} + e_i) \\
&=& \frac{1}{N^d}  \int_{I_N} A_{per} \left(\nabla w_i^{0} + e_i \right) \\
&=& A_{per}^{*} e_i,
\end{eqnarray*}
which obviously gives the zero-order term expected for $A_{\eta}^*$. Then
\begin{eqnarray}\label{a1starn}
 A_{1}^{*,N} e_i &=& \frac{1}{N^d}   \sum_{l=1}^{N^d} \left \langle d\bar{P}_1(s_l)\prod_{k=1, k\neq l}^{N^d} \delta_0(s_k), \int_{I_N} A^{s_1, \cdots, s_{N^d}} ( \nabla w_i^{s_1, \cdots, s_{N^d}} + e_i ) \right \rangle.
\end{eqnarray}
It is easy to see that, by $(N \ZZ)^d$-periodicity of $w_i^{s_1, \cdots, s_{N^d}}$,
$$\left \langle d\bar{P}_1(s_l)\prod_{k=1, k\neq l}^{N^d} \delta_0(s_k), \int_{I_N} A^{s_1, \cdots, s_{N^d}} ( \nabla w_i^{s_1, \cdots, s_{N^d}} + e_i ) \right \rangle$$
does not depend on $l$. The expression (\ref{a1starn}) can then be rewritten
\begin{eqnarray}\label{order1aux}
 A_{1}^{*,N} e_i &=& \left \langle  d\bar{P}_1(s), \int_{I_N} A^{s, 0 \cdots, 0} ( \nabla w_i^{s,0, \cdots, 0} + e_i ) \right \rangle.  
\end{eqnarray}

We change the notations for convenience, and define, for $s \in [-M,M]$,
\begin{eqnarray} \label{a1s0}
 A_1^{s,0} = A^{s, 0 \cdots, 0} = A_{per} + s \mathds{1}_Q C_{per},
\end{eqnarray}
 and $w_i^{1,s,0,N} = w_i^{s,0, \cdots, 0}$ solution to 
\begin{eqnarray}\label{sdefect} 
\left  \{
\begin{aligned}
& -\mathrm{div}\left(A_1^{s,0}(\nabla w_i^{1,s,0,N} + e_i ) \right) = 0 \quad \, \mathrm{in}  \, I_N,\\
& w_i^{1,s,0,N} \,(N\ZZ)^d-\mathrm{periodic}.
\end{aligned}
\right.
\end{eqnarray} 

The matrix $A_1^{s,0}$ corresponds to the periodic material with a defect of amplitude $s$ located in $Q$ (i.e at a position $0 \in \ZZ^d$ in $I_N$), and $w_i^{1,s,0,N}$ is the $i$-th cell solution for the periodic homogenization of $A_1^{s,0}$ in $I_N$. Since $w_i^{1,s,0,N} = w_i^{s,0, \cdots, 0}$, it is of course a $C^{\infty}$ function of $s \in ]-M,M[$.\\

With these notations, we find that
\begin{eqnarray} \label{order1}
  A_{1}^{*,N} e_i = \left \langle  d\bar{P}_1(s), \int_{I_N} A_1^{s,0} (\nabla w_i^{1,s,0,N} + e_i )  \right \rangle.
\end{eqnarray}

For the second-order term, we first define the set
\begin{eqnarray} \label{taun}
\mathcal{T}_N = \left\{ k \in \ZZ^d, Q+k \subset I_N \right\} = \left \llbracket -\frac{N-1}{2}, \frac{N-1}{2} \right \rrbracket^d.
\end{eqnarray}
The cardinal of $\mathcal{T}_N$ is of course $N^d$, and $\displaystyle \bigcup_{k \in \mathcal{T}_N} \{Q+k\} = I_N.$ 

For $(s,t) \in [-M,M]^2$ and $k \in \mathcal{T}_N$, we define
\begin{eqnarray}
A_2^{s,t,0,k} = A_{per} + s \mathds{1}_Q C_{per} + t \mathds{1}_{Q + k} C_{per},
\end{eqnarray}
and $w_i^{2,s,t,0,k,N}$ solution to 
\begin{eqnarray} \label{sdefect2}
\left  \{
\begin{aligned}
& -\mathrm{div}\left(A_2^{s,t,0,k}(\nabla w_i^{2,s,t,0,k,N} + e_i ) \right) = 0 \quad \, \mathrm{in}  \, I_N,\\
& w_i^{2,s,t,0,k,N} \,(N\ZZ)^d-\mathrm{periodic}.
\end{aligned}
\right.
\end{eqnarray}

The matrix $A_2^{s,t,0,k}$ corresponds to the periodic material with two defects of amplitude $s$ and $t$ located in $Q$ and $Q+k$ (i.e at positions $0 \in \ZZ^d$ and $k \in \ZZ^d$ in $I_N$) respectively. The function $w_i^{2,s,t,0,k,N}$ is the $i$-th cell solution for the periodic homogenization of $A_2^{s,t,0,k}$ in $I_N$. It is a $\mathcal{C}^{\infty}$ function of $(s,t) \in ]-M,M[^2$.\\

Then computations similar to that presented for the first order yield
\begin{eqnarray} \label{order2}
\begin{aligned}
A_{2}^{*,N} e_i = \frac{1}{2} \sum_{k \in \mathcal{T}_N, k \neq 0}\left \langle d\bar{P}_1(s) d\bar{P}_1(t) ,  \int_{I_N} A_2^{s,t,0,k}(\nabla w_i^{2,s,t,0,k,N} + e_i) \  \right \rangle \\
 + \left \langle d\bar{P}_2(s), \int_{I_N}A_1^{s,0} (\nabla w_i^{1,s,0,N} + e_i) \right \rangle. 
\end{aligned}
\end{eqnarray}

\begin{figure}[H]
\center
\begin{tabular}{ccc}
\includegraphics[width=4.5cm, height=4.5cm]{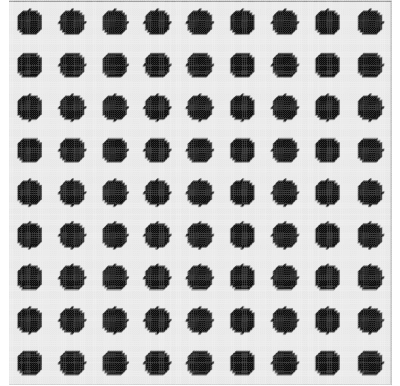} & \includegraphics[width=4.5cm, height=4.5cm]{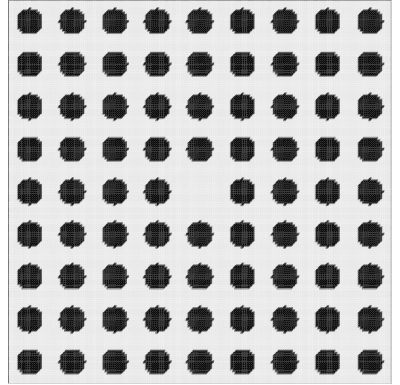} & \includegraphics[width=4.5cm, height=4.5cm]{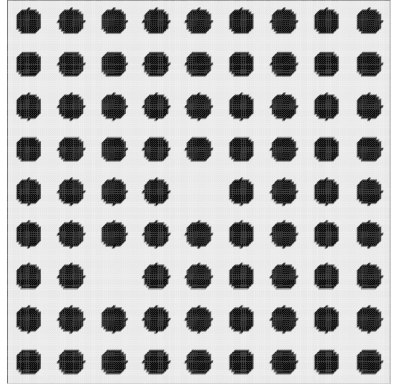}
\end{tabular}
\caption{From left to right: zero defect, one defect and two defects.}
\label{illus}
\end{figure}

A setting with zero, one and two defects is shown in Figure \ref{illus} in the two-dimensional case of a reference material $A_{per}$ consisting of a periodic lattice of circular inclusions.

\begin{remark}
It is illustrative to consider the particular case where the random variable $B_{\eta}$ has a Bernoulli law. This is the case treated in \cite{ALB1}. Then, expansion (\ref{pushexpand}) holds exactly with $d\bar{P}_1 = \delta_1 - \delta_0$. The distribution $d\bar{P}_2$ and all other terms of higher order identically vanish. The expressions (\ref{order1}) and (\ref{order2}) then coincide with (3.17) and (3.18) in \cite{ALB1}.
\end{remark}

In the next section we prove that $A_{1}^{*,N}$ converges to a finite limit when $N \rightarrow \infty$. The case of the second-order term  $A_{2}^{*,N}$, which is shown to be a bounded sequence and thus to converge up to extraction, is discussed in Section \ref{seco}.

\subsection{Convergence of the first-order term}

We study here the convergence as $N$ goes to infinity of $A_{1}^{*,N}$ defined by (\ref{order1}). 

\begin{Prop} \label{conva1}
The sequence $A_{1}^{*,N}$ converges in $\mathbb{R}^{d \times d}$ to a finite limit $\bar{A_1^*}$ when $N \rightarrow \infty$.
\end{Prop}

\begin{proof} We fix $(i,j) \in \llbracket 1,d \rrbracket^2$ and study the convergence of $A_{1}^{*,N} e_i \cdot e_j$.\\

Using (\ref{sdefect}) and the adjoint problems defined by (\ref{celladj}), we first obtain, for all \newline $s\in [-M,M]$,
\begin{eqnarray*}
  \int_{I_N} A_1^{s,0} (\nabla w_i^{1,s,0,N} + e_i) \cdot e_j  =\int_{I_N} A_1^{s,0} (\nabla w_i^{1,s,0,N} + e_i ) \cdot (e_j + \nabla \tilde{w}_{j}^0). 
\end{eqnarray*}
Then, letting the distribution $d\bar{P}_1$ act on the left and right-hand sides, and using (\ref{order1}), we find that
\begin{eqnarray} \label{order1mod}
  A_{1}^{*,N} e_i \cdot e_j = \left \langle  d\bar{P}_1(s), \int_{I_N} A_1^{s,0} (\nabla w_i^{1,s,0,N} + e_i) \cdot (e_j + \nabla \tilde{w}_{j}^0)\right \rangle.
\end{eqnarray}

Because of the definition of $A_1^{s,0}$, 
\begin{eqnarray} \label{step1}
\begin{aligned}
\int_{I_N} A_1^{s,0} (\nabla w_i^{1,s,0,N} + e_i) \cdot (e_j + \nabla \tilde{w}_{j}^0) = \int_{I_N} A_{per} (\nabla w_i^{1,s,0,N} + e_i) \cdot (e_j + \nabla \tilde{w}_{j}^0) \\ + \int_{Q} sC_{per}(\nabla w_i^{1,s,0,N} + e_i) \cdot (e_j + \nabla \tilde{w}_{j}^0).
\end{aligned}
\end{eqnarray}

Next, using (\ref{celladj}),
\begin{eqnarray} \label{step2}
\begin{aligned}
\int_{I_N}A_{per}(\nabla w_{i}^{1,s,0,N} + e_i ) \cdot (e_j + \nabla \tilde{w}_{j}^0) &=  \int_{I_N}(\nabla w_{i}^{1,s,0,N} + e_i ) \cdot A_{per}^T (e_j + \nabla \tilde{w}_{j}^0) \\
&= \int_{I_N} e_i \cdot A_{per}^T (e_j + \nabla \tilde{w}_{j}^0). 
\end{aligned}
\end{eqnarray}

We know from Lemma \ref{distreasy} that $\langle d\bar{P}_1, 1\rangle=0$. Thus
\begin{eqnarray} \label{step3}
 \left \langle  d\bar{P}_1(s), \int_{I_N} e_i \cdot A_{per}^T (e_j + \nabla \tilde{w}_{j}^0)\right \rangle = 0.
\end{eqnarray}

Collecting (\ref{order1mod}), (\ref{step1}), (\ref{step2}) and (\ref{step3}), we get
\begin{eqnarray} \label{adef1genaux}
A_{1}^{*,N} e_i \cdot e_j = \left \langle  d\bar{P}_1(s), \int_{Q} sC_{per}(\nabla w_i^{1,s,0,N} + e_i) \cdot (e_j + \nabla \tilde{w}_{j}^0) \right \rangle.
\end{eqnarray}

We now define
\begin{eqnarray} \label{diff1}
q_i^{1,s,0,N} =w_i^{1,s,0,N} -w_i^{0}.
\end{eqnarray}

$q_i^{1,s,0,N}$ solves
\begin{eqnarray} \label{perturbdefgen} 
\left  \{
\begin{aligned}
& -\mathrm{div} \left(A_1^{s,0} \nabla q_i^{1,s,0,N} \right) = \mathrm{div} (s \mathds{1}_Q C_{per} (\nabla w_i^0 + e_i)) \quad \mathrm{in} \; I_N,\\
& q_i^{1,s,0,N} \;  (N\ZZ)^d-\mathrm{periodic}.
\end{aligned}
\right.
\end{eqnarray}

Using (\ref{diff1}) in (\ref{adef1genaux}), we rewrite
\begin{eqnarray} \label{adef1gen}
A_{1}^{*,N} e_i \cdot e_j &=& \left \langle  s d\bar{P}_1(s), \int_{Q}C_{per}(\nabla w_i^{0} + e_i) \cdot (e_j + \nabla \tilde{w}_{j}^0) \right \rangle  \nonumber \\
&& \qquad + \left \langle  d\bar{P}_1(s), \int_{Q} sC_{per}(\nabla q_i^{1,s,0,N} + e_i) \cdot (e_j + \nabla \tilde{w}_{j}^0) \right \rangle.
\end{eqnarray}

The rest of the proof consists in showing that 
$$ \left \langle  d\bar{P}_1(s), \int_{Q} sC_{per}(\nabla q_i^{1,s,0,N} + e_i ) \cdot (e_j + \nabla \tilde{w}_{j}^0) \right \rangle,$$
which is of course equal to
 $$\left \langle   s d\bar{P}_1(s), \int_{Q} C_{per}(\nabla q_i^{1,s,0,N} + e_i) \cdot (e_j + \nabla \tilde{w}_{j}^0) \right \rangle,$$
converges to a finite limit when $N \rightarrow \infty$.\\

More precisely, defining
$$\forall s \in [-M,M], \; \forall N \in 2 \NN +1 , \; \; f^N(s) = \int_{Q} C_{per}(\nabla q_i^{1,s,0,N} + e_i) \cdot (e_j + \nabla \tilde{w}_{j}^0),$$
we will prove that the sequence $f^N$ and its derivatives converge uniformly, when $N$ goes to infinity, to a limit function $f^{\infty}$ and its derivatives.\\

Applying Lemma \ref{conv1} of the appendix to (\ref{perturbdefgen}), we obtain that for all $s \in [-M,M]$, 
$\nabla q_i^{1,s,0,N}$ converges in $L^2(Q)$, when $N \rightarrow \infty$, to $\nabla q_i^{1,s,0,\infty}$, where $q_i^{1,s,0,\infty}$ is a $L^2_{loc}(\RR^d)$ function solving
\begin{eqnarray} \label{perturbdefgeninf} 
\left  \{
\begin{aligned}
& -\mathrm{div} \left(A_1^{s,0} \nabla q_i^{1,s,0,\infty} \right) = \mathrm{div} (s \mathds{1}_Q C_{per} (\nabla w_i^0 + e_i)) \quad \mathrm{in} \; \RR^d,\\
& \nabla q_i^{1,s,0,\infty} \in L^2(\RR^d).
\end{aligned}
\right.
\end{eqnarray}

Moreover, arguing as in the proof of Lemma \ref{conv1} (given in our previous work \cite{ALB1}), it is easy to see that for all $k \in \NN$ and all $s \in [-M,M]$, $\nabla \partial_s^k q_i^{1,s,0,N}$ converges in $L^2(Q)$ to  $\nabla \partial_s^k q_i^{1,s,0,\infty}$.\\

We then define $f^\infty$ by
$$\forall s \in [-M,M], \; \; f^\infty(s) = \int_{Q} C_{per}(\nabla q_i^{1,s,0,\infty} + e_i) \cdot (e_j + \nabla \tilde{w}_{j}^0).$$

Because of (\ref{bornder}) and (\ref{bornderinf}) in Lemma \ref{distribq1} of the appendix, and using a classical result of differentiation under the integral sign, it is clear that
$$ \forall k \in \NN, \forall s \in ]-M,M[, \; \frac{d^k}{ds^k} f^N(s) = \int_{Q} C_{per}(\nabla \partial^k_s q_i^{1,s,0,N} + e_i) \cdot (e_j + \nabla \tilde{w}_{j}^0),$$
and
$$ \forall k \in \NN, \forall s \in ]-M,M[, \; \frac{d^k}{ds^k} f^{\infty}(s) = \int_{Q} C_{per}(\nabla \partial^k_s q_i^{1,s,0,\infty} + e_i) \cdot (e_j + \nabla \tilde{w}_{j}^0).$$

The convergence of $\nabla \partial_s^k q_i^{1,s,0,N}$ to  $\nabla \partial_s^k q_i^{1,s,0,\infty}$ in $L^2(Q)$ for every $k \in \NN$ thus yields
\begin{eqnarray} \label{convsimple}
\forall k \in \NN, \forall s \in ]-M,M[, \lim_{N \rightarrow + \infty }\frac{d^k}{ds^k} f^N(s)  =\frac{d^k}{ds^k} f^\infty(s). 
\end{eqnarray}

On the other hand, we deduce from Lemma \ref{uniflipschitz} that there exists a constant $C(p_1,M)$ (recall that $p_1$ is the order of $d\bar{P}_1(s)$) such that
for all $k \in \llbracket 0, p_1 \rrbracket$,
\begin{eqnarray} \label{lipschitzo}
 \forall (s,s') \in ]-M,M[^2, \forall N \in 2 \NN+1, \; \; |\frac{d^k}{ds^k} f^N(s) - \frac{d^k}{ds^k} f^N(s')| \leq C(p_1,M) |s-s'|. 
\end{eqnarray}

It is straightforward to see that (\ref{convsimple}) and (\ref{lipschitzo}) imply that
\begin{eqnarray} \label{convuniforme}
 \forall 0 \leq k \leq p_1, \;\frac{d^k}{ds^k} f^N \mathrm{\;converges \; uniformly \; to \;} \frac{d^k}{ds^k} f^\infty \; \mathrm{\;in \;} ]-M,M[.
\end{eqnarray}

It follows from (\ref{major1}) and (\ref{convuniforme}) that 
$$ \langle s d\bar{P}_1(s) , f^N(s) \rangle \rightarrow \langle s d\bar{P}_1(s) , f^{\infty}(s) \rangle,$$
and then
\begin{eqnarray}
\begin{aligned} \label{convdp1}
& \left \langle  d\bar{P}_1(s), \int_{Q} sC_{per}(\nabla q_i^{1,s,0,N} + e_i) \cdot (e_j + \nabla \tilde{w}_{j}^0) \right \rangle\\
 & \qquad \qquad \underset{N \rightarrow \infty}{\rightarrow}  \left \langle  d\bar{P}_1(s), \int_{Q} sC_{per}(\nabla q_i^{1,s,0,\infty} + e_i) \cdot (e_j + \nabla \tilde{w}_{j}^0)\right \rangle.
\end{aligned}
\end{eqnarray}

Collecting (\ref{adef1gen}) and (\ref{convdp1}), we conclude that $A_1^{*,N}$ converges to a limit tensor $\bar{A_1^*}$ defined by
\begin{eqnarray} \label{ainfiny}
\begin{aligned}
\forall (i,j) \in \llbracket 1, d \rrbracket^2, \; \; \bar{A_1^*} e_i \cdot e_j &= \left \langle  s d\bar{P}_1(s), \int_{Q}C_{per}\left(\nabla w_i^{0} + e_i \right) \cdot (e_j + \nabla \tilde{w}_{j}^0) \right \rangle  \\
&   + \left \langle  d\bar{P}_1(s), \int_{Q} sC_{per}\left(\nabla q_i^{1,s,0,\infty} + e_i \right) \cdot (e_j + \nabla \tilde{w}_{j}^0) \right \rangle.
\end{aligned}
\end{eqnarray}

\end{proof}

\subsection{Second-order term} \label{seco}

For completeness, we state here the result concerning the second-order term $A_{2}^{*,N}$ in (\ref{asymptN}), proved in \cite{these_Arnaud}:
\begin{Prop} \label{conva2}
The sequence $A_{2}^{*,N}$ defined by (\ref{order2}) is bounded in $\RR^{d\times d}$ and therefore converges up to extraction.
\end{Prop}

We firmly believe that $A_{2}^{*,N}$ is actually a convergent sequence, as shown by our numerical tests thereafter. We also stress that the explicit computations of Section \ref{onedimension} prove this convergence in dimension one.

\section{Numerical experiments}\label{num}

The purpose of this section is to assess the numerical relevance of the approaches of Sections \ref{model} and \ref{heuristic}. To this end we build and homogenize stochastic composite materials using laws that satisfy the assumptions of these sections. Our motivations are not strictly identical for the two approaches. In contrast to the first approach which relies on a rigorous proof, our second approach is formal and we thus need to demonstrate its correctness experimentally (note that the tests performed in \cite{ALB1} in the Bernoulli case are already to be considered as a component of the validation of the approach). We wish to check that the expansions derived in Sections \ref{model} and \ref{heuristic} provide an accurate and efficient approximation to the direct stochastic computation. The limited computational facilities we have access to impose that we restrict ourselves to the two-dimensional setting. We first explain our general methodology, which is the same as that presented in \cite{ALB1}, and then make precise the specific settings.

\subsection{Methodology}

We mainly consider as in \cite{ALB1} a reference material $A_{per}$ that consists of a constant background reinforced by a periodic lattice of circular inclusions, that is
$$ A_{per}(x_1,x_2) = 20 \times Id + 100 \sum_{k \in \ZZ^2} \mathds{1}_{B(k,0.3)}(x_1,x_2) \times Id,$$
where $B(k,0.3)$ is the ball of center $k$ and radius $1$. Loosely speaking, the role of the perturbation is to randomly eliminate some fibers:
$$ C_{per}(x_1,x_2) = -100 \sum_{k \in \ZZ^d} \mathds{1}_{B(k,0.3)}(x_1,x_2) \times Id.$$

We will also, in our last test, consider a laminate $$ A_{per}(x_1,x_2) = 5 + 10 \sum_{l \in \ZZ} \mathds{1}_{l \leq x_1 \leq l+1}(x_1,x_2),$$  with the perturbation yielding an error in the lamination direction:
$$ C_{per}(x_1,x_2) =  10 \sum_{l \in \ZZ} \mathds{1}_{l \leq x_2 \leq l+1}(x_1,x_2)\times Id- 10 \sum_{l \in \ZZ} \mathds{1}_{l \leq x_1 \leq l+1}(x_1,x_2) \times Id.$$

\begin{figure}[h] 
\center
\begin{tabular}{cc}
\includegraphics[width = 7cm, height = 7cm]{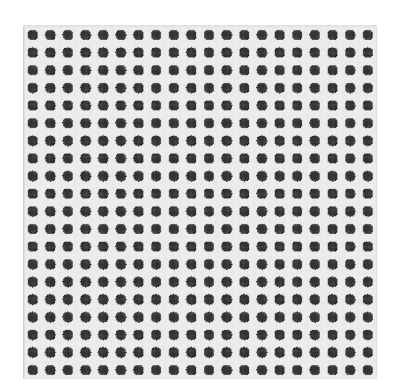} & \includegraphics[width = 7cm, height = 7cm]{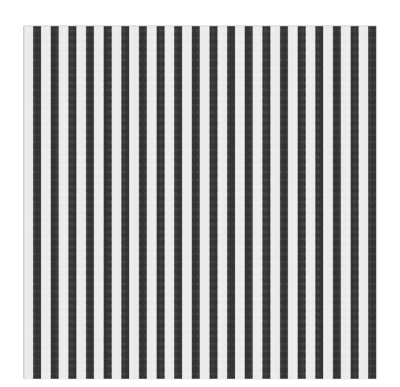}
\end{tabular}
\caption{Left: a periodic lattice of circular inclusions. Right: a one-dimensional laminate.}
\label{material}
\end{figure}

For both materials (shown in Figure \ref{material}), we have chosen the values of the coefficients in order to have a high contrast between $A_{per}$ and $A_{per} + C_{per}$ and thus for the perturbation to have an important impact on the microscopic structure. The specific value of these coefficients has no other significance.  \\

We will consider different perturbations $b_{\eta}$, all of which satisfy (\ref{forme}) with the $B_{\eta}^k$ independent and identically distributed.\\

Our goal is to compare $A_{\eta}^*$ with its approximation $A_{per}^{*} + \eta A_1^{*,N} + \eta^2 A_2^{*, N}$. A major computational difficulty is the computation of the ``exact'' matrix $A_{\eta}^*$ given by formula (\ref{homogal}). It ideally requires to solve the stochastic cell problems (\ref{cellsol}) on $\RR^d$. To this end we first use ergodicity and formula (\ref{expect}), and actually compute, for a given realization $\omega$ and a domain $I_N$ chosen here to be $[0,N]^2$ for convenience, $A_{\eta}^{*,N}(\omega)$ defined by 
\begin{eqnarray} \label{ergodicint} 
A_{\eta}^{*,N}(\omega) e_i = \frac{1}{N^d} \int_{I_N}A_{\eta}(x,\omega)(\nabla w_{i}^{\eta,N,\omega}(x) + e_i ) dx. 
\end{eqnarray} 

In a second step, we take averages over the realizations $\omega$.\\

For each $\omega$, we use the finite element software FreeFem++ (available at www.freefem.org) to solve the boundary value problems (\ref{ergodiccorr}) and compute the integrals (\ref{ergodicint}). We work with standard P1 finite elements on a triangular mesh such that there are $10$ degrees of freedom on each edge of the unit cell $Q$.\\  

We define an approximate value $A_{\eta}^{*,N}$ as the average of $A_{\eta}^{*,N}(\omega)$ over $40$ realizations~$\omega$. Our numerical experiments indeed show that the number $40$ is sufficiently large for the convergence of the Monte-Carlo computation. We then let $N$ grow from $5$ to $80$ by steps of $5$. We observe that $A_{\eta}^{*,N}$ stabilizes at a fixed value around $N=80$ and thus take $A_{\eta}^{*,80}$ as the reference value for $A_{\eta}^{*}$ in our subsequent tests.\\

The next step is to compute the zero-order term $A_{per}^*$, and the first-order and second-order deterministic corrections. Using the same mesh and finite elements as for our reference computation above, we compute $A_{per}^*$ using (\ref{cellper}) and (\ref{homper}). The computation of the next orders depends on the setting:
\begin{itemize}
\item in the setting of Section $2$, the first-order correction is given by (\ref{tilda1}) in Theorem \ref{rig} and is thus independent of $N$; since $b_{\eta}$ is of the form (\ref{specific}), we use formula (\ref{expans2correl}) in Corollary \ref{correl} for the second-order correction which depends on $N$ through the term $t_i$ defined on $\RR^d$ by (\ref{ti}), and which has to be approximated on $I_N$; we let $N$ grow from $5$ to $80$ by steps of $5$;
\item in the setting of Section $3$, the corrections $A_1^{*,N}$ and $A_2^{*,N}$ are respectively given by (\ref{order1}) and (\ref{order2}); we let $N$ grow from $5$ to $80$ by steps of $5$ for $A_1^{*,N}$; the computation of $A_2^{*,N}$ being far more expensive (there is not only an integral over $I_N$ but also a sum over the $N^2$ cells in (\ref{order2})), we have to limit ourselves to $N=25$ and approximate the value for $N$ larger than $25$ by the value obtained for $N=25$.
\end{itemize}

We stress that there are three distinct sources of error in these computations:
\begin{itemize}
\item the finite elements discretization error; 
\item the truncation error due to the replacement of $\RR^d$ with $I_N$, in the computation of the stochastic cell problems (\ref{cellsol}) that are replaced with (\ref{ergodiccorr}), as well as in the computation of the integrals (\ref{ergodicint}); 
\item the stochastic error arising from the approximation of the expectation value by an empirical mean. 
\end{itemize}

Detailed comments on these various errors and the way we deal with them are provided in \cite{ALB1}. We just emphasize, in the setting of Section \ref{heuristic}, that it is not our purpose to prove through our tests that    
$$ A_{\eta}^* = A_{per}^* + \eta \bar{A_1^*} + \eta^2 \bar{A_2^*} + o(\eta^2)$$
with a $o(\eta^2)$ which would be independent of $N$, of the number of realizations and of the size of the mesh. We only wish to demonstrate that the second-order expansion is an approximation to $A_{\eta}^*$ sufficiently good for all practical purposes. We will observe that $A_2^{*,N}$ is not only bounded as stated in Proposition \ref{conva2} but converges to a limit $\bar{A_2^*}$, and that both $A_1^{*,N}$ and $A_2^{*,N}$ converge to their respective limits faster than $A_{\eta}^{*,N}$ to $A_{\eta}^{*}$ (which is expected since the former quantities are deterministic and contain less information). We will also observe that $A_{per}^* + \eta A_1^{*,N}$ is closer to $A_{\eta}^{*}$ than $A_{per}^{*}$ and that the inclusion of the second order improves the situation for $A_{per}^* + \eta A_1^{*,N} + \eta^2 A_2^{*,N}$ is even closer.\\

To present our numerical results, we choose the first diagonal entry $(1,1)$ of all the matrices considered. Other coefficients in the matrices behave qualitatively similarly. We illustrate a practical interval of confidence for our Monte-Carlo computation of $A_{\eta}^*$ by showing, for each $N$, the minimum and maximum values of $A_{\eta}^{*,N}(\omega)$ achieved over the $40$ realizations $\omega$.\\

We will use the following legend in the graphs:
\begin{itemize}
\item {\it periodic}: gives the value of the periodic homogenized tensor $A_{per}^*$;
\item {\it first-order}: gives the value of the first-order expansion;
\item {\it second-order}: gives the value of the second-order expansion;
\item {\it stochastic mean, minima and maxima}: respectively give the values of $A_{\eta}^{*,N}$ and the extrema obtained in the computation of the empirical mean.\\
\end{itemize}

Finally, the results are given for various values of $\eta$ which serve the purpose of testing our approach in a diversity of situations, and in particular for perturbations that are ``not so small''. 

\subsection{An example of setting for our theory in Section \ref{model} (and \ref{heuristic})}

Consider $B_{\eta} = \eta \, G \, \mathds{1}_{0 \leq \eta G  \leq 1}$ where $G$ is a normalized centered Gaussian random variable. It is easy to check that
$$B_{\eta} = \eta G \mathds{1}_{0 \leq G  \leq + \infty} + o(\eta^2) \quad \mathrm{in \;} L^2(\Omega),$$ 
so that Corollary \ref{correl} of Section 2 applies. Alternatively, we can use Lemma \ref{lemexp},  which gives
\begin{eqnarray*}
dP_{\eta} &=& \delta_0 - \eta \frac{1}{\sqrt{2}}\delta_0'+ \frac{\eta^2}{4}\delta_0'' + o(\eta^2) \; \mathrm{in} \; \mathcal{E}'(\mathbb{R}),
\end{eqnarray*}
to perform our formal approach. We verify in Section \ref{consistency} of the appendix that both approaches yield the same results up to second order.\\

\begin{figure}
\center
\includegraphics[width=14cm, height=11cm]{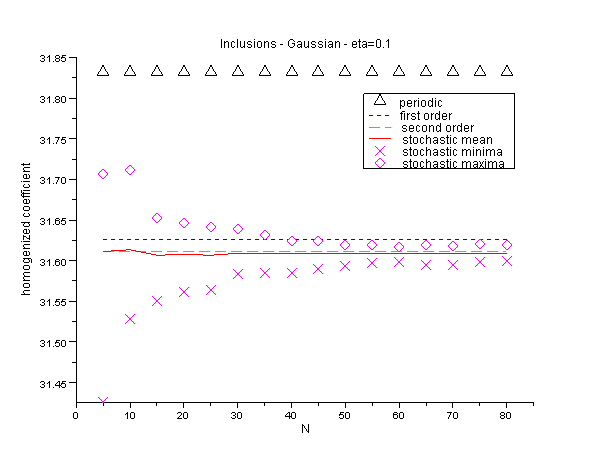} 
\includegraphics[width=14cm, height=11cm]{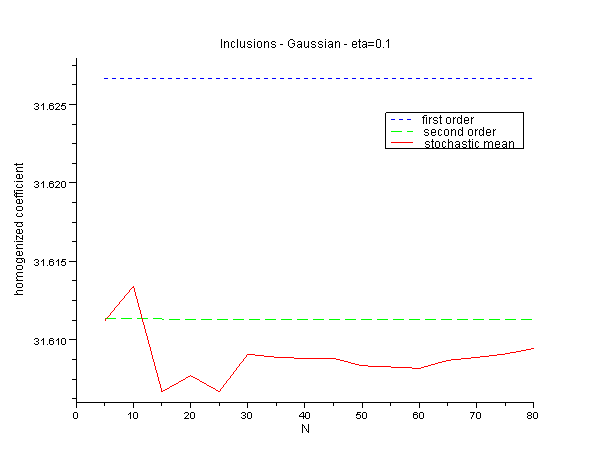} 
\caption{Inclusions - results for a Gaussian perturbation and $\eta=0.1$. Above: complete results. Below: close-up on $A_{\eta}^{*,N}$ and the first and second-order corrections.}
\label{convfibres_G_01}
\end{figure}

\begin{figure}
\center
\includegraphics[width=14cm, height=11cm]{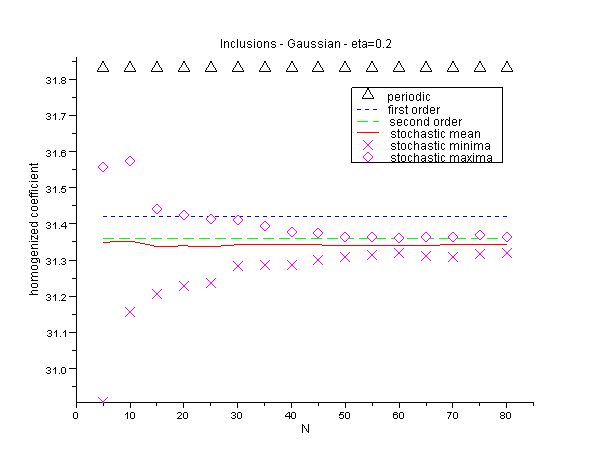} 
\includegraphics[width=14cm, height=11cm]{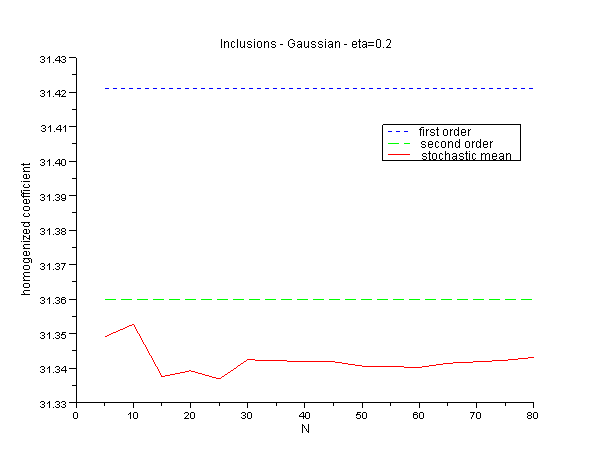} 
\caption{Inclusions - Results for a Gaussian perturbation and $\eta=0.2$. Above: complete results. Below: close-up on $A_{\eta}^{*,N}$ and the first and second-order corrections.}
\label{convfibres_G_02}
\end{figure}

We show results for the lattice of inclusions and for $\eta=0.1$ and $\eta=0.2$ (Figures \ref{convfibres_G_01} and \ref{convfibres_G_02} respectively).\\

The results are very satisfying for both values of $\eta$. The first-order correction, which does not depend on $N$, enables to get substantially closer to $A_{\eta}^*$. Moreover, it is clear (especially from the close-ups) that the second-order correction $A_2^{*,N}$ converges very fast (convergence is already reached at $N=5$), and in particular much faster than the stochastic computation $A_{\eta}^{*,N}$. It also provides excellent accuracy.

\subsection{A first example of setting for our formal approach of Section \ref{heuristic}}

Consider $R_{\eta}$ a random variable having Bernoulli law with parameter $\eta$, and $G$ a normalized centered Gaussian random variable independent of $R_{\eta}$. We define the product random variable $B_{\eta} = R_{\eta} \times \eta G \mathds{1}_{|\eta G|\leq 1}$. Then
\begin{eqnarray*}
\mathbb{E}(\varphi(B_{\eta})) &=& \mathbb{E}(\varphi(R_{\eta} \times \eta G \mathds{1}_{|\eta G|\leq 1})) \\ 
&=& \eta \mathbb{E}(\varphi(\eta G \mathds{1}_{|\eta G|\leq 1})) + (1-\eta) \varphi(0) \\
&=& \eta (\varphi(0) + \eta \mathbb{E}(G) \varphi'(0) + \frac{\eta^2}{2} \varphi''(0) + o(\eta^2)) + (1-\eta) \varphi(0) \\
&=& \varphi(0) + \frac{\eta^3}{2} \varphi''(0) + o(\eta^3).
\end{eqnarray*}
 
This implies 
\begin{eqnarray} \label{order3}
dP_{\eta} = \delta_0 + \frac{\eta^3}{2} \delta^{''}_0 + o(\eta^3) \quad \mathrm{in \;} \mathcal{E}'(\RR).
\end{eqnarray}

\begin{figure}
\center
\includegraphics[width=14cm, height=11cm]{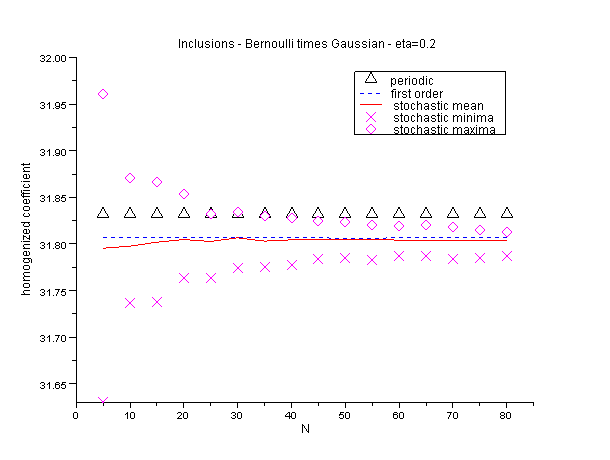} 
\includegraphics[width=14cm, height=11cm]{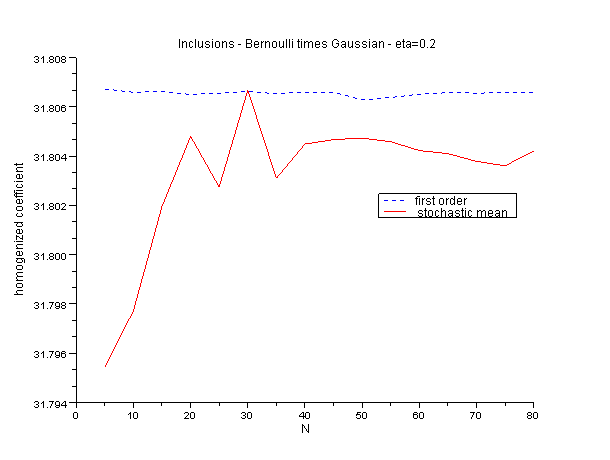} 
\caption{Inclusions - results for perturbation (\ref{order3}) and $\eta=0.1$. Above: complete results. Below: close-up on $A_{\eta}^{*,N}$ and the first-order correction.}
\label{convfibres_BG_02}
\end{figure}

\begin{figure}
\center
\includegraphics[width=14cm, height=11cm]{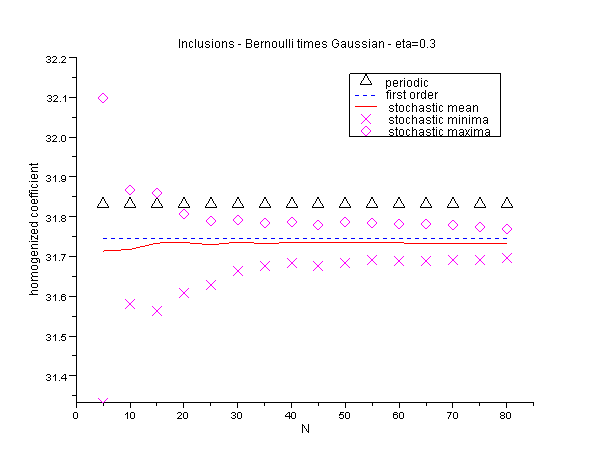} 
\includegraphics[width=14cm, height=11cm]{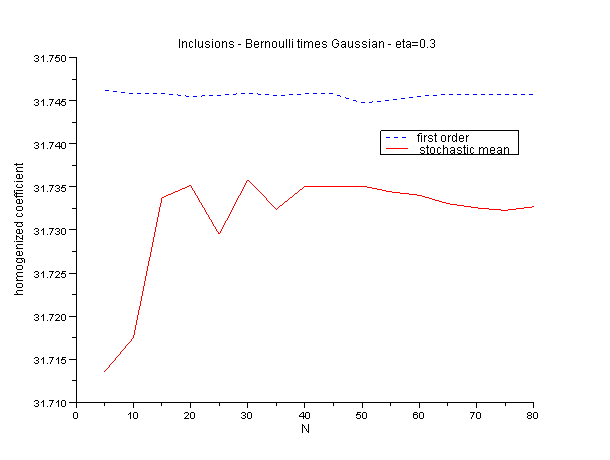} 
\caption{Inclusions - results for perturbation (\ref{order3}) and $\eta=0.3$. Above: complete results. Below: close-up on $A_{\eta}^{*,N}$ and the first-order correction.}
\label{convfibres_BG_03}
\end{figure}

\begin{figure}
\center
\includegraphics[width=14cm, height=11cm]{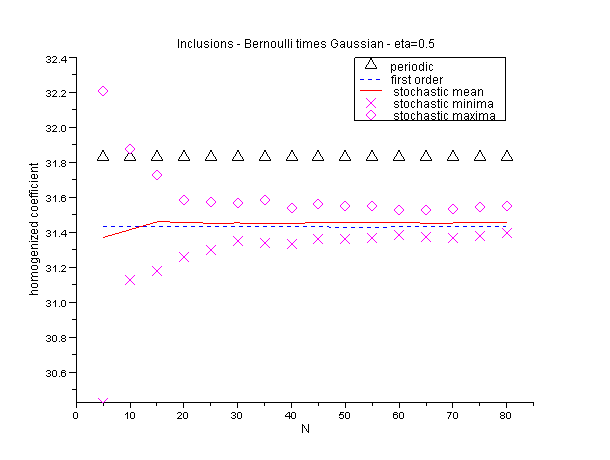} 
\includegraphics[width=14cm, height=11cm]{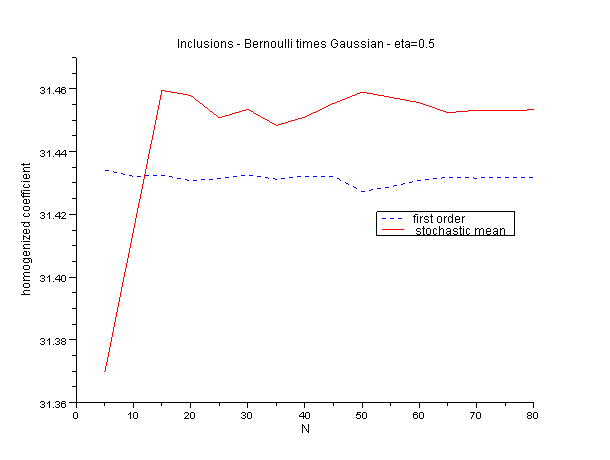} 
\caption{Inclusions - results for perturbation (\ref{order3}) and $\eta=0.5$. Above: complete results. Below: close-up on $A_{\eta}^{*,N}$ and the first-order correction.}
\label{convfibres_BG_05}
\end{figure}

In this case we only consider the first-order correction since the dominant order in (\ref{order3}) is already tiny. We present the results in the case of the lattice of inclusions, for $\eta=0.2$, $\eta=0.3$ and $\eta=0.5$ (Figures \ref{convfibres_BG_02}, \ref{convfibres_BG_03}, \ref{convfibres_BG_05} respectively).\\

Once again, our approach converges rapidly and allows for an accurate approximate value of $A_{\eta}^{*}$ even for $\eta$ as large as 0.5. 

\newpage

\subsection{A second example of setting for our formal approach of Section \ref{heuristic}}

Consider $R_{\eta}$ a random variable having Bernoulli law with parameter $\eta$, and $U$ a uniform variable on $[0,1]$ independent of $R_{\eta}$.  We define $B_{\eta} = R_{\eta} - \eta U$. Then
\begin{eqnarray*}
\mathbb{E}(\varphi(B_{\eta})) &=& \mathbb{E}(\varphi(R_{\eta} - \eta U)) \\
&=& \eta \mathbb{E}(\varphi(1 - \eta U)) + (1-\eta) \mathbb{E}(\varphi(-\eta U))  \\
&=& \eta \left(\varphi(1) - \eta \mathbb{E}(U) \varphi'(1) + o(\eta)\right) \\
 &&+ (1-\eta) \left(\varphi(0) - \eta \mathbb{E}(U) \varphi'(0)  + \frac{\eta^2}{2}\mathbb{E}(U^2) \varphi''(0) + o(\eta^2) \right)\\
&=& \varphi(0) + \eta \left(-\mathbb{E}(U) \varphi'(0) + \varphi(1)-\varphi(0)\right) \\
&&+ \eta^2\left(-\mathbb{E}(U)(\varphi'(1)-\varphi'(0)) + \frac{1}{2}\mathbb{E}(U^2) \varphi''(0)\right) + o(\eta^2),
\end{eqnarray*}
so that 
\begin{eqnarray} \label{complex}
\begin{aligned}
dP_{\eta} =& \delta_0 + \eta \left(-\mathbb{E}(U) \delta'_0 + \delta_1-\delta_0 \right)\\
 &+ \eta^2\left(-\mathbb{E}(U)(\delta'_1-\delta'_0) + \frac{1}{2}\mathbb{E}(U^2) \delta''(0)\right) + o(\eta^2) \quad \mathrm{in \;} \mathcal{E}'(\RR).
\end{aligned}
\end{eqnarray}

Notice that this complex case is a mixture of Sections 2 and 3. The first-order perturbation is of course only the sum of the first-order perturbations for a Bernoulli law (Section 3 and \cite{ALB1}) and a uniform law (Section 2). The interaction of these laws at order 2, and notably the $\delta'_1$ term, is much more involved and requires the computation of the cross derivatives of $w_i^{2,s,t,0,k,N}$ with respect to $s$ and $t$ at $s=0$ and $t=1$. \\

\begin{figure}[H]
\center
\includegraphics[width=14cm, height=11cm]{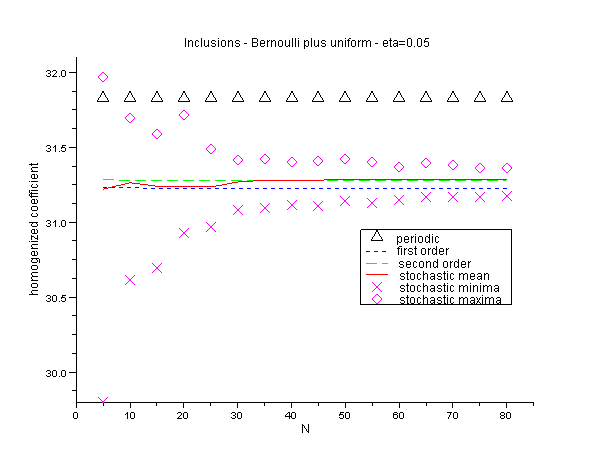} 
\includegraphics[width=14cm, height=11cm]{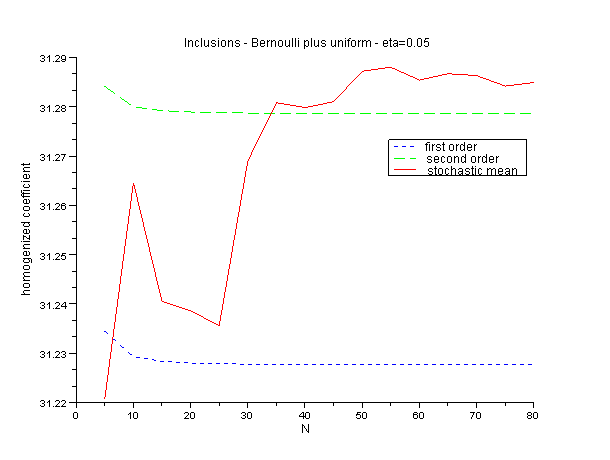} 
\caption{Inclusions- results for perturbation (\ref{complex}) and $\eta=0.05$. Above: complete results. Below: close-up on $A_{\eta}^{*,N}$ and the first and second-order corrections.}
\label{convfibres_BU_005}
\end{figure}

\begin{figure}[H]
\center
\includegraphics[width=14cm, height=11cm]{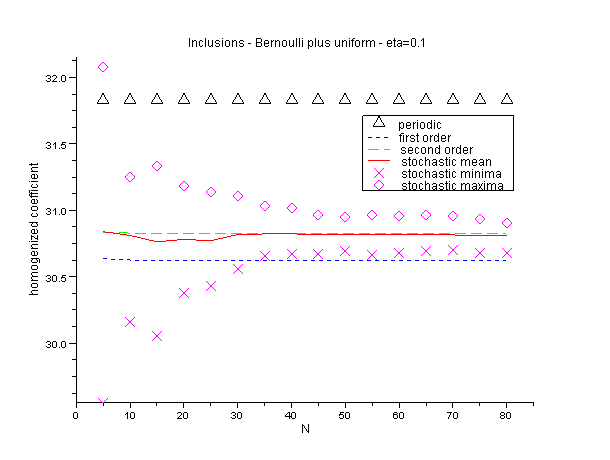} 
\includegraphics[width=14cm, height=11cm]{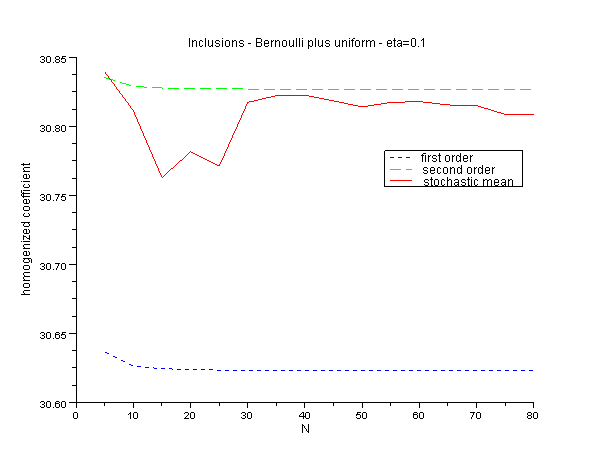} 
\caption{Inclusions - results for perturbation (\ref{complex}) and $\eta=0.1$. Above: complete results. Below: close-up on $A_{\eta}^{*,N}$ and the first and second-order corrections.}
\label{convfibres_BU_01}
\end{figure}

\begin{figure}[H]
\center
\includegraphics[width=14cm, height=11cm]{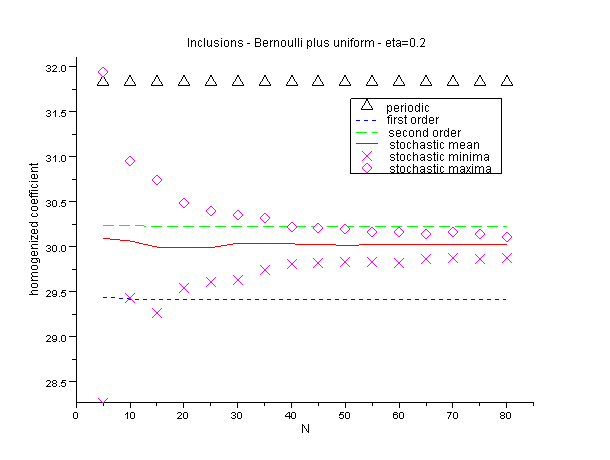} 
\includegraphics[width=14cm, height=11cm]{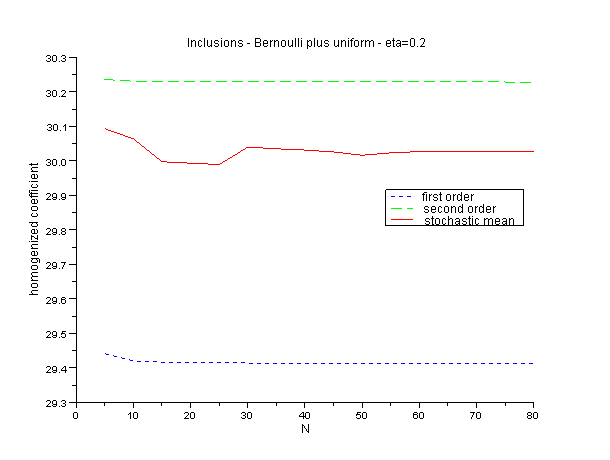} 
\caption{Inclusions - results for perturbation (\ref{complex}) and $\eta=0.2$. Above: complete results. Below: close-up on $A_{\eta}^{*,N}$ and the first and second-order corrections.}
\label{convfibres_BU_02}
\end{figure}

We give the results in the case of the inclusions and for $\eta=0.05$, $\eta=0.1$ and $\eta=0.2$ (Figures \ref{convfibres_BU_005}, \ref{convfibres_BU_01}, \ref{convfibres_BU_02}, respectively). \\

For $\eta=0.05$ and $\eta=0.1$, the results display the same features as in our previous tests and are very good. The case $\eta=0.2$ is instructive: the second-order expansion significantly departs from the "exact" value provided by the direct stochastic computation. Our interpretation is that, far from contradicting the validity of our expansion in the limit of small $\eta$, it shows the limitations of the approach. The value $\eta=0.2$ is too large for the expansion to be accurate in the case of a lattice of inclusions with a high contrast between the inclusions and the surrounding phase.\\

\begin{figure}[H]
\center
\includegraphics[width=14cm, height=11cm]{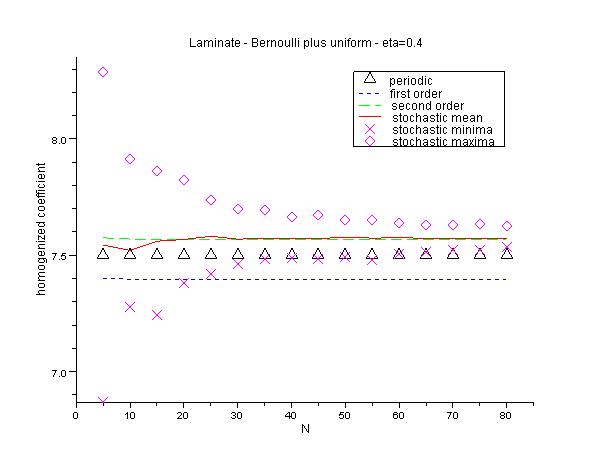} 
\includegraphics[width=14cm, height=11cm]{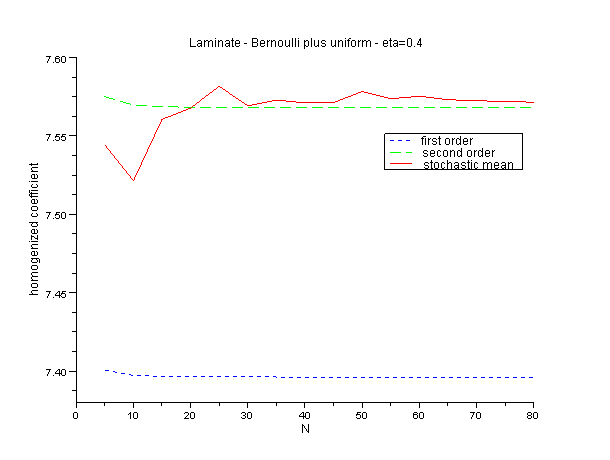} 
\caption{Laminate - results for perturbation (\ref{complex}) and $\eta=0.4$. Above: complete results. Below: close-up on $A_{\eta}^{*,N}$ and the first and second-order corrections.}
\label{convlamine_BU_04}
\end{figure}

Interestingly, a value of $\eta$ twice as large (0.4) provides a very accurate approximation for another material, as shown by our final test performed on the laminate (Figure \ref{convlamine_BU_04}).\\

 Our approach has limitations and deteriorates, like any asymptotic approach, for large values of $\eta$. The threshold is case dependent. The approach is however generically robust.

\newpage

\newpage

\section{Appendix}

The objectives of this appendix are diverse. We first quickly recall some elements of distribution theory. We then prove technical results used in Section \ref{heuristic}. Next we show that the approach formally derived in Section \ref{heuristic} is rigorous in dimension one. Finally we prove that this approach is also rigorous, in general dimensions, in a specific setting close to that of Theorem \ref{rig} and Corollary \ref{correl}.

\subsection{Elements of distribution theory} \label{distheo}

We recall here some basic definitions and results of distribution theory for convenience of the reader. See \cite{Hörmander} for a comprehensive presentation.\\

In this section $\mathcal{O}$ denotes an open set in $\RR$.

\begin{defi}
We denote by $\mathcal{D}(\mathcal{O})$ the space of infinitely differentiable functions on $\mathcal{O}$ having compact support in $\mathcal{O}$. 
\end{defi}

\begin{defi}
T is a distribution on $\mathcal{O}$ if $T$ is a linear form on $\mathcal{D}(\mathcal{O})$ satisfying the following continuity property: for every compact $K \subset \mathcal{O}$, there exists an integer $p$ and a constant $C$ such that for all $\varphi \in \mathcal{D}(\mathcal{O})$ having compact support in $K$,
\begin{eqnarray} \label{defidist}
|\langle T, \varphi \rangle| \leq C \sup_{x \in K, 0 \leq n \leq p} \left|\frac{d^n}{dx^n} \varphi(x) \right|.
\end{eqnarray}
The space of distributions on $\mathcal{O}$ is denoted by $\mathcal{D}'(\mathcal{O})$.\\

If the integer $p$ in (\ref{defidist}) can be chosen independently of $K$, the distribution $T$ is said to have a finite order. The smallest possible value for $p$ is called the order of $T$.  
\end{defi}

\begin{defi}
A distribution $T \in \mathcal{D}'(\mathcal{O})$ is said to have compact support if there exists a compact set $K \subset \mathcal{O}$ such that for all $\varphi \in \mathcal{D}(\mathcal{O})$ having compact support in $\mathcal{O} \backslash K$, $\langle T, \varphi \rangle = 0$.\\

The support of $T$ is defined as the smallest compact set $K$ which satisfies the above assertion.\\

The space of distributions on $\mathcal{O}$ having compact support is denoted by $\mathcal{E}'(\mathcal{O})$. 

\end{defi}

\begin{Prop} \label{extendcont}
If $T \in \mathcal{E}'(\mathcal{O})$, its action on $\mathcal{D}(\mathcal{O})$ can be naturally extended to $\mathcal{C}^{\infty}(\mathcal{O})$. Denoting by $K$ a compact neighborhood of the support of $T$, and by $\chi$ a cut-off function in $\mathcal{D}(\mathcal{O})$ equal to 1 on the support of $T$ and vanishing on $\mathcal{O} \backslash K$, we define
$$ \forall \varphi \in \mathcal{C}^{\infty}(\mathcal{O}), \; \langle T, \varphi \rangle := \langle T, \chi \varphi \rangle.$$

This definition does not depend on $K$ and $\chi$.
\end{Prop}

\begin{Prop} \label{finiteorder}
If a distribution $T$ is in $\mathcal{E}'(\mathcal{O})$, it has a finite order. Denoting by $p$ its order and by $K$ a compact neighborhood of the support of $T$, there exists a constant $C>0$ such that:
$$ \forall \varphi \in C^{\infty}(\mathcal{O}), \; |\langle T, \varphi \rangle| \leq C \sup_{x \in K, 0 \leq n \leq p} \left|\frac{d^n}{dx^n} \varphi(x) \right|.$$
\end{Prop}


\subsection{Some technical results} \label{technical}

This section is devoted to the proof of technical lemmas used in Section \ref{heuristic}. Loosely speaking, these lemmas all deal with the variation of the supercell correctors defined by (\ref{sgen}), (\ref{sdefect}), and (\ref{sdefect2}) with respect to the amplitudes of the defects.

\begin{lemma} \label{regular}
Let $\tilde{H}^1_{per}(I_N)$ be the set of $(N \ZZ)^d$-periodic functions in $H^1_{loc}(\RR^d)$ with zero mean on $I_N$. 
The function $$F : (s_1, \cdots, s_{N^d}) \in ]-M,M[^{N^d} \mapsto  \bar{w}_i^{s_1, \cdots, s_{N^d}} \in \tilde{H}^1_{per}(I_N),$$
where $\bar{w}_i^{s_1, \cdots, s_{N^d}} = w_i^{s_1, \cdots, s_{N^d}} - \int_{I_N} w_i^{s_1, \cdots, s_{N^d}}$ and $w_i^{s_1, \cdots, s_{N^d}}$ is defined by (\ref{sgen}),  is $C^{\infty}$.
\end{lemma}

\begin{proof}
For $(s_1, \cdots, s_{N^d}) \in [-M,M]^{N^d}$, $\bar{w}_i^{s_1, \cdots, s_{N^d}}$ is the unique solution to
\begin{eqnarray*} 
\left  \{
\begin{aligned}
& -\mathrm{div}\left(A^{s_1, \cdots, s_{N^d}}(\nabla \bar{w}_i^{s_1, \cdots, s_{N^d}}+ e_i ) \right) = 0 \quad \, \mathrm{in}  \, I_N,\\
& \bar{w}_i^{s_1, \cdots, s_{N^d}} \,(N\ZZ)^d-\mathrm{periodic}, \quad \int_{I_N} \bar{w}_i^{s_1, \cdots, s_{N^d}} = 0,
\end{aligned}
\right.
\end{eqnarray*} 

so that $F$ is well defined.

Let us now define $G : ]-M,M[^{N^d} \times \tilde{H}^1_{per}(I_N) \rightarrow H^{-1}(I_N)$ by
$$ G(s_1, \cdots, s_{N^d},w) =  -\mathrm{div}\left(A^{s_1, \cdots, s_{N^d}}(\nabla w+ e_i ) \right),$$
so that
$F(s_1, \cdots, s_{N^d}) = \bar{w}_i^{s_1, \cdots, s_{N^d}}$ is the unique solution to 
$$ G(s_1, \cdots, s_{N^d}, F(s_1, \cdots, s_{N^d}))=0.$$

It is easy to see that $G$ is a $C^1$ function, and that
$$\forall h \in \tilde{H}^1_{per}(I_N), \quad \partial_w G(s_1, \cdots, s_{N^d},w) \cdot h = -\mathrm{div}\left(A^{s_1, \cdots, s_{N^d}}\nabla h \right),$$
where $\partial_w G(s_1, \cdots, s_{N^d},w)$ is the first derivative of $G$ with respect to $w$ at $(s_1, \cdots, s_{N^d},w)$.\\

The Lax-Milgram theorem and the coercivity of $A^{s_1, \cdots, s_{N^d}}$ show that $\partial_w G(s_1, \cdots, s_{N^d},w)$ is an isomorphism. We can therefore apply the inverse function theorem and deduce that $F$ is $C^1$, with $\partial_{s_l} F$ the unique solution to 
\begin{eqnarray*} 
\left  \{
\begin{aligned}
& -\mathrm{div}\left(A^{s_1, \cdots, s_{N^d}}(\nabla \partial_{s_l} F ) \right) = \mathrm{div} \left(\mathds{1}_{Q_l} C_{per} (\nabla F + e_i)    \right) \quad \, \mathrm{in}  \, I_N,\\
& \partial_{s_l} F \,(N\ZZ)^d-\mathrm{periodic}, \quad \int_{I_N} \partial_{s_l} F = 0.
\end{aligned}
\right.
\end{eqnarray*} 

Arguing by induction, we obtain that $F$ is a $C^{\infty}$ function.

\end{proof}

For consistency, we state next a lemma proved in \cite[Lemma 6]{ALB1} 

\begin{lemma} \label{conv1}
Consider $f \in L^2(Q)$, and a tensor field $A$ from $\RR^d$ to $\RR^{d\times d}$ such that there exist $\lambda>0$ and $\Lambda > 0$ such that
$$ \forall \xi \in \RR^d, \mathrm{\; a.e \; in \;} x \in \RR^d, \; \lambda |\xi|^2 \leq A(x)\xi\cdot \xi \; \mathrm{and} \; |A(x) \xi| \leq \Lambda |\xi|.$$ 
 Consider $q^N$ solution to 
\begin{eqnarray} \label{perturblem}
\left  \{
\begin{aligned}
& -\mathrm{div} \left(A \nabla q^{N} \right) = \mathrm{div} (\mathds{1}_Q f) \quad \mathrm{in} \; I_N,\\
& q^{N}\,(N \ZZ)^d-\mathrm{periodic}.
\end{aligned}
\right.
\end{eqnarray} 
Then $ \mathds{1}_{I_N} \nabla q^N$ converges in $L^2(\RR^d)$, when $N$ goes to infinity, to $\nabla q^{\infty}$, where $q^{\infty}$ is a~$L^2_{loc}(\RR^d)$ function solving
\begin{eqnarray} \label{perturblemdef}  
\left  \{
\begin{aligned}
& -\mathrm{div} \left(A \nabla q^{\infty} \right) = \mathrm{div} (\mathds{1}_Q f) \quad \mathrm{in} \; \RR^d,\\
& \nabla q^{\infty} \in L^2(\RR^d).
\end{aligned}
\right.
\end{eqnarray} 
\end{lemma}

\begin{lemma} \label{distribq1}
Consider $q_i^{1,s,0,N}$ and $q_i^{1,s,0,\infty}$ solutions to (\ref{perturbdefgen}) and (\ref{perturbdefgeninf}) respectively, and $k \in \NN$. There exists a constant $C(k,M)$, such that
\begin{eqnarray} \label{bornder}
\forall s \in ]-M,M[\,, \; \forall N \in 2 \NN + 1, \; \; \|\nabla \partial_s^{k} q_i^{1,s,0,N}\|_{L^2(I_N)} \leq C(k,M) \|\nabla w_i^{0} + e_i\|_{{L^2}(Q)},
\end{eqnarray}
\begin{eqnarray} \label{bornderinf}
\forall s \in ]-M,M[\,, \; \; \|\nabla \partial_s^{k} q_i^{1,s,0,\infty}\|_{L^2(\RR^d)} \leq C(k,M) \|\nabla w_i^{0} + e_i\|_{{L^2}(Q)}.
\end{eqnarray}
\end{lemma}

\begin{proof}
Multiplying the first line of (\ref{perturbdefgen}) by $q_i^{1,s,0,N}$ and integrating by parts, we find that
\begin{eqnarray} \label{bornder0}
\|\nabla q_i^{1,s,0,N}\|_{L^2(I_N)} \leq M \frac{\|C_{per}\|_{L^{\infty}(Q)}}{\alpha} \|\nabla w_i^{0} + e_i\|_{{L^2}(Q)},
\end{eqnarray}
where $\alpha$ is defined by (\ref{alphabound}).\\

Thus (\ref{bornder}) is true for $k=0$ with $C(0,M) = M \frac{\|C_{per}\|_{L^{\infty}(Q)}}{\alpha}$.\\

Next, the first derivative $\partial_s q_i^{1,s,0,N}$ is solution to 
\begin{eqnarray} \label{der1q}
\left  \{
\begin{aligned}
& -\mathrm{div}(A_1^{s,0} \nabla \partial_s q_i^{1,s,0,N}) = \mathrm{div}\left(\mathds{1}_{Q}C_{per}(\nabla w_i^{0} + e_i)\right) + \mathrm{div}\left(\mathds{1}_{Q}C_{per}\nabla q_i^{1,s,0,N}\right) \quad \mathrm{in} \; I_N,\\
& \partial_s q_i^{1,s,0,N} \; \; (N\ZZ)^d-\mathrm{periodic},
\end{aligned}
\right.
\end{eqnarray}

from which we deduce
\begin{eqnarray*}
\|\nabla \partial_s q_i^{1,s,0,N}\|_{L^2(I_N)} \leq \frac{\|C_{per}\|_{L^{\infty}(Q)}}{\alpha} \left(\|\nabla w_i^{0} + e_i\|_{{L^2}(Q)}+\|\nabla q_i^{1,s,0,N}\|_{{L^2}(Q)} \right)
\end{eqnarray*}
and, using (\ref{bornder0}),
\begin{eqnarray*}
\|\nabla \partial_s q_i^{1,s,0,N}\|_{L^2(I_N)} \leq \frac{\|C_{per}\|_{L^{\infty}(Q)}}{\alpha} (M+1) \|\nabla w_i^{0} + e_i\|_{{L^2}(Q)}.
\end{eqnarray*}

Thus (\ref{bornder}) is true for $k=1$ with $C(1,M) = (M+1) \frac{\|C_{per}\|_{L^{\infty}(Q)}}{\alpha}$.\\

Finally, we have for $k \geq 2$
\begin{eqnarray} \label{der1qk}
\left  \{
\begin{aligned}
& -\mathrm{div}(A_1^{s,0} \nabla \partial_s^k q_i^{1,s,0,N}) = k \mathrm{div}\left(\mathds{1}_{Q}C_{per}\nabla \partial_s^{k-1} q_i^{1,s,0,N}\right) \quad \mathrm{in} \; I_N,\\
& \nabla \partial_s^k q_i^{1,s,0,N} \; \; (N\ZZ)^d-\mathrm{periodic},
\end{aligned}
\right.
\end{eqnarray}

so that an easy induction proves (\ref{bornder}). The proof of (\ref{bornderinf}) is identical.
\end{proof}

The following result is an immediate consequence of Lemma \ref{distribq1}.

\begin{lemma}\label{uniflipschitz}
Consider $q_i^{1,s,0,N}$ and $q_i^{1,s,0,\infty}$ solutions to (\ref{perturbdefgen}) and (\ref{perturbdefgeninf}) respectively. For every $k \in \NN$, there exists a constant $C(k,M)$ such that for all $(s,s') \in ]-M,M[^2$,
\begin{eqnarray} \label{lipschitzlem}
\forall N \in 2 \NN + 1, \; \; \|\nabla \partial_s^k q_i^{1,s,0,N} - \nabla \partial_s^k q_i^{1,s',0,N}\|_{L^2(I_N)} \leq C(k,M)\|\nabla w_i^{0} + e_i\|_{{L^2}(Q)} |s-s'|, 
\end{eqnarray}
\begin{eqnarray} \label{lipschitzleminf}
 \|\nabla \partial_s^k q_i^{1,s,0,\infty} - \nabla \partial_s^k q_i^{1,s',0,\infty}\|_{L^2(\RR^d)} \leq C(k,M)\|\nabla w_i^{0} + e_i\|_{{L^2}(Q)} |s-s'|. 
\end{eqnarray}
\end{lemma}

%

\subsection{The one-dimensional case} \label{onedimension}

We address here the one-dimensional context. All the computations are explicit, for the settings of Sections 2 and \ref{heuristic}. To stress the fact that we deal with scalar quantities, we use lower-case letters for the tensors. Note also that in this section $Q=[-\frac{1}{2},\frac{1}{2}]$ and $I_N = [-\frac{N}{2},\frac{N}{2}]$.

\subsubsection{An extension of Theorem \ref{rig}}

The following theorem extends the result of Theorem \ref{rig}, stated in $L^{\infty}(Q;L^2(\Omega))$, \newline to $L^{\infty}(Q;L^p(\Omega))$ for any $p \in ]1, \infty]$:
\begin{theorem}[one-dimensional setting] \label{rig1d}
Assume that $d=1$, that $b_{\eta}$ satisfies (\ref{hyp1}) and $m_{\eta} := \|b_{\eta}\|_{L^{\infty}([-\frac{1}{2},\frac{1}{2}]; L^{p}(\Omega))} \underset{\eta \rightarrow 0}{\rightarrow 0}$ for some $p>1$. There exists a subsequence of $\eta$, still denoted $\eta$ for simplicity, such that $\frac{b_{\eta}}{m_{\eta}}$ converges weakly-* in ${L^{\infty}([-\frac{1}{2},\frac{1}{2}]; L^{p}(\Omega))}$ to a limit field denoted by $\bar{b}_0$ when $\eta \rightarrow 0$. Then 
\begin{itemize}
\item the expansion
\begin{eqnarray} \label{expancell1d}
\frac{d}{dx}w^{\eta} = \frac{d}{dx}w^0 + m_{\eta} \frac{d}{dx} v^0 + o(m_{\eta})
\end{eqnarray}
holds weakly in $L^2([-\frac{1}{2},\frac{1}{2}];L^p(\Omega))$, 
where $w^0$ is the periodic corrector and $v^{0}$ solves 
\begin{eqnarray} \label{vound}
\left  \{
\begin{aligned}
& -\frac{d}{dx}(a_{per} \frac{d}{dx} v^{0}) = \frac{d}{dx}\left(\bar{b}_0c_{per}(\frac{d}{dx} w^{0} + 1)\right) \quad \mathrm{in} \; \RR,\\
& \frac{d}{dx} v^{0} \; \mathrm{stationary}, \; \mathbb{E}\left(\int_{-\frac{1}{2}}^{\frac{1}{2}} \frac{d}{dx} v^{0}\right) = 0.
\end{aligned}
\right.
\end{eqnarray}
\item $a_{\eta}^*$ reads 
\begin{eqnarray*} \label{expans1d}
a_{\eta}^{*}  = a_{per}^{*}+ m_{\eta}\int_{-\frac{1}{2}}^{\frac{1}{2}}  \mathbb{E}(\bar{b}_0)c_{per} (\frac{d}{dx} w^{0} +  1) + m_{\eta}\int_{-\frac{1}{2}}^{\frac{1}{2}} a_{per}\frac{d}{dx}  \mathbb{E}(v^{0}) + o(m_{\eta}).
\end{eqnarray*}
\end{itemize}
\end{theorem}

\begin{proof}

The periodic and stochastic correctors can be computed explicitly. They are respectively given by
$$\frac{d}{dx} w^0 = \left(\int_{-\frac{1}{2}}^{\frac{1}{2}} a_{per}^{-1} \right)^{-1} a_{per}^{-1} -1 \; \; \mathrm{and} \; \; \frac{d}{dx} w^{\eta} = \left(\mathbb{E}\left(\int_{-\frac{1}{2}}^{\frac{1}{2}} a_{\eta}^{-1}\right) \right)^{-1} a_{\eta}^{-1} -1.$$

Note that $w^0$ is in $W^{1,\infty}(-\frac{1}{2},\frac{1}{2})$.\\

We define $\displaystyle v^{\eta} = \frac{w^{\eta} - w^0}{m_{\eta}}$. It solves
\begin{eqnarray} \label{vund}
\left  \{
\begin{aligned}
& -\frac{d}{dx}(a_{\eta} \frac{d}{dx} v^{\eta}) = \frac{d}{dx}\left(\frac{b_{\eta}}{\eta} c_{per}(\frac{d}{dx} w^{0} + 1) \right)  \quad \mathrm{in} \; \RR,\\
& \frac{d}{dx} v^{\eta} \; \; \mathrm{stationary}, \; \; \mathbb{E}\left(\int_{-\frac{1}{2}}^{\frac{1}{2}} \frac{d}{dx} v^{\eta}\right) = 0.
\end{aligned}
\right.
\end{eqnarray}
We deduce from (\ref{vund}) that
\begin{eqnarray} \label{vund2}
a_{\eta} \frac{d}{dx} v^{\eta} = \frac{b_{\eta}}{m_{\eta}}c_{per}(\frac{d}{dx} w^{0} + 1) + k_{\eta},
\end{eqnarray}
where $k_{\eta}$ depends only on $\omega$. Since $k_{\eta}$ is by construction stationary ergodic, it is constant, and we compute from (\ref{vund}) and (\ref{vund2}):
$$k_{\eta} = - \frac{1}{m_{\eta}}\left(\mathbb{E}\int_{-\frac{1}{2}}^{\frac{1}{2}} \frac{1}{a_{\eta}}\right)^{-1} \times \left(\mathbb{E}\int_{-\frac{1}{2}}^{\frac{1}{2}} \frac{b_{\eta}}{a_{\eta}} c_{per}(\frac{d}{dx} w^{0} + 1) \right).$$

Since $w^0$ is in $W^{1,\infty}(-\frac{1}{2},\frac{1}{2})$, $a_{\eta}$ is coercive and $c_{per}$ is bounded, it holds
\begin{eqnarray*}
\left|k_{\eta}\right| &\leq& C \frac{\|b_{\eta}\|_{L^1([-\frac{1}{2},\frac{1}{2}] \times \Omega)}}{m_{\eta}} \\
&\leq& C \frac{\|b_{\eta}\|_{L^1([-\frac{1}{2},\frac{1}{2}] \times \Omega)}}{\|b_{\eta}\|_{L^{\infty}([-\frac{1}{2},\frac{1}{2}] ; L^p(\Omega))}}.
\end{eqnarray*}

This implies that $k_{\eta}$ is a bounded function of $\eta$ whatever $p\geq1$ and thus, using (\ref{vund2}), that $\frac{d}{dx} v^{\eta}$ is bounded in $L^2([-\frac{1}{2},\frac{1}{2}]; L^p(\Omega))$ for all $p \geq 1$. As a result, for $p>1$,  $\frac{d}{dx} v^{\eta}$ converges weakly and up to extraction in $L^2([-\frac{1}{2},\frac{1}{2}]; L^p(\Omega))$ to a limit we denote $\frac{d}{dx} v^0$.\\

The random field $b_{\eta}$ tends to $0$ in $L^2([-\frac{1}{2},\frac{1}{2}]; L^p(\Omega))$. Since it is bounded in $L^{\infty}([-\frac{1}{2},\frac{1}{2}]\times \Omega)$, it converges to $0$ in $L^2([-\frac{1}{2},\frac{1}{2}]; L^r(\Omega))$ for all $r > p$. By Hölder inequality it also converges to $0$ in $L^2([-\frac{1}{2},\frac{1}{2}]; L^r(\Omega))$ for all $1 < r <p$. Thus it converges to $0$ in in $L^2([-\frac{1}{2},\frac{1}{2}]; L^q(\Omega))$ where $q = \frac{p}{p-1}$.\\

The space $L^2([-\frac{1}{2},\frac{1}{2}]; L^q(\Omega))$ being the dual of $L^2([-\frac{1}{2},\frac{1}{2}]; L^p(\Omega))$, we obtain that $b_{\eta} \, c_{per} \, \frac{d}{dx} v^{\eta}$ tends to $0$ in $\mathcal{D}'([-\frac{1}{2},\frac{1}{2}] \times \Omega)$. We can then take the limit $\eta \rightarrow 0$ in (\ref{vund}) and obtain that $v^0$ is solution to (\ref{vound}).\\

We have thus proved that $\frac{1}{m_{\eta}} \left(\frac{d}{dx} w^{\eta} - \frac{d}{dx} w^0\right)$ converges, up to extraction, weakly to $\frac{d}{dx} v^0$ in $L^2([-\frac{1}{2},\frac{1}{2}];L^p(\Omega))$, which is equivalent to (\ref{expancell1d}).\\

The second assertion of Theorem \ref{rig1d} is obtained by inserting (\ref{expancell1d}) into the expression (\ref{homogal}) of $a_{\eta}^*$.
\end{proof}

Note that the proof of Theorem \ref{rig1d} depends crucially on the fact that we are able to solve explicitly the cell problems.\\

Theorem \ref{rig1d} allows for a better intuitive understanding of Theorem \ref{rig}. In dimension one, the homogenized coefficient is explicitly given by
$$a_{\eta}^{*} = \left(\mathbb{E}\int_{-\frac{1}{2}}^{\frac{1}{2}} \frac{1}{a_{per} + b_{\eta} c_{per}}\right)^{-1},$$
which, when $ \displaystyle b_{\eta}(x,\omega) = \sum_{k \in \ZZ} \mathds{1}_{[k,k+1]}(x) B_{\eta}(\tau_k \omega)$, may be rewritten as the \emph{formal} series
\begin{eqnarray}\label{aux}
\frac{1}{a_{\eta}^{*}} = \sum_{k=0}^{\infty} (-1)^k \mathbb{E}((B_{\eta})^k)  \int_{-\frac{1}{2}}^{\frac{1}{2}} \left(\frac{c_{per}}{a_{per}}\right)^k a_{per}^{-1}.
\end{eqnarray}

Assume now that there exists $p>1$ such that $\|B_{\eta}\|_{L^{p}(\Omega)} \rightarrow 0 $ when $\eta \rightarrow 0$ and $\frac{B_{\eta}}{\|B_{\eta}\|_{L^{p}(\Omega)}}$ converges weakly in $L^p(\Omega)$ to some $\bar{B}_0$ with $\mathbb{E}(\bar{B}_0) \neq 0$. We have in particular
 $$ \frac{\mathbb{E}(B_{\eta})} {\|B_{\eta}\|_{L^{p}(\Omega)}}  \rightarrow \mathbb{E}(\bar{B}_0)\neq 0,$$
which, since $\mathbb{E}(|B_{\eta}|^p) \rightarrow 0$, implies 
\begin{eqnarray} \label{petito}
\mathbb{E}(|B_{\eta}|^p) = \underset{\eta \rightarrow 0^+}{o}\left( \mathbb{E}(B_{\eta}) \right).
\end{eqnarray}
We now claim that, without loss of generality and up to an extraction in $\eta$, we may take $p=2$ in (\ref{petito}). Indeed, if $p<2$, then since $B_{\eta}$ is bounded in $L^{\infty}(\Omega)$, (\ref{petito}) implies $\mathbb{E}(|B_{\eta}|^2) = \underset{\eta \rightarrow 0^+}{o}\left( \mathbb{E}(B_{\eta}) \right)$. On the other hand, if $p>2$, we consider the normalized sequence $\frac{B_{\eta}}{\|B_{\eta}\|_{L^{2}(\Omega)}}$ in $L^2(\Omega)$. Up to extraction, it weakly converges to $\bar{B}_2 \in L^2(\Omega)$. Since
$$  \frac{\mathbb{E}(B_{\eta})} {\|B_{\eta}\|_{L^{p}(\Omega))}}  = \frac{\mathbb{E}(B_{\eta})} {\|B_{\eta}\|_{L^{2}(\Omega)}} \frac{\|B_{\eta}\|_{L^{2}(\Omega)}}{\|B_{\eta}\|_{L^{p}(\Omega)}}$$
where the left hand side converges to  $\mathbb{E}(\bar{B}_0)\neq 0$ and $\frac{\|B_{\eta}\|_{L^{2}(\Omega)}}{\|B_{\eta}\|_{L^{p}(\Omega)}}$ is bounded by 1 by Hölder's inequality, $ \displaystyle \mathbb{E}(\bar{B}_2) = \lim_{\eta \rightarrow 0}\frac{\mathbb{E}(B_{\eta})} {\|B_{\eta}\|_{L^{2}(\Omega))}}\neq 0$ and (\ref{petito}) is satisfied with $p=2$.\\

We then take $p=2$. Since $\mathbb{E}(|B_{\eta}|^2) = \underset{\eta \rightarrow 0^+}{o}\left( \mathbb{E}(B_{\eta}) \right)$ and $B_{\eta}$ is bounded in $L^{\infty}(\Omega)$, $\mathbb{E}(|B_{\eta}|^k) = \underset{\eta \rightarrow 0^+}{o}\left( \mathbb{E}(B_{\eta}) \right)$ for all $k \geq 2$.\\

This intuitively expresses that all orders higher than or equal to 2 are negligible as compared to the first-order term in the series (\ref{aux}), and thus that a kind of ``separation of scales'' is satisfied. This is of course formal since one has to check that the remainder term consisting of the \emph{sum} of all terms of order higher than or equal to 2 is $o\left( \mathbb{E}(B_{\eta}) \right)$, so that 
\begin{eqnarray*}
a_{\eta}^{*} &=& \left(\int_{-\frac{1}{2}}^{\frac{1}{2}} \frac{1}{a_{per}}\right)^{-1} +\left(\int_{-\frac{1}{2}}^{\frac{1}{2}} \frac{1}{a_{per}}\right)^{-2} \left(\int_{-\frac{1}{2}}^{\frac{1}{2}} \mathbb{E}(B_{\eta}) \frac{c_{per}}{a_{per}}\right) + o\left( \mathbb{E}(B_{\eta}) \right)\\
&=& \left(\int_{-\frac{1}{2}}^{\frac{1}{2}} \frac{1}{a_{per}}\right)^{-1} +m_{\eta}\left(\int_{-\frac{1}{2}}^{\frac{1}{2}} \frac{1}{a_{per}}\right)^{-2} \left(\int_{-\frac{1}{2}}^{\frac{1}{2}} \mathbb{E}(\bar{B}_0) \frac{c_{per}}{a_{per}}\right) + o\left( \mathbb{E}(B_{\eta}) \right).
\end{eqnarray*}
But this is the purpose of the proofs of Theorems \ref{rig} and \ref{rig1d}, using another viewpoint, to show this is indeed the case.

\subsubsection{The setting of Section \ref{heuristic} in dimension one} 

We now prove that our approach of Section \ref{heuristic} is rigorous in dimension one. 

\begin{lemma}
In dimension $d=1$, it holds
$$a_{\eta}^{*} = a_{per}^* + \eta \bar{a}_1^* + \eta^2 \bar{a}_2^* + o(\eta^2),$$
where $\bar{a}_1^*$ and $\bar{a}_2^*$ are the limits as $N \rightarrow \infty$ of $a_1^{*,N}$ and $a_2^{*,N}$ defined generally by (\ref{order1}) and (\ref{order2}) respectively. 
\end{lemma}

\begin{proof}

Recall that in dimension one, $a_{\eta}^*$ is given by the simple explicit expression
$$a_{\eta}^{*} = \left(\mathbb{E}\int_{-\frac{1}{2}}^{\frac{1}{2}} \frac{1}{a_{per} + b_{\eta} c_{per}}\right)^{-1} =  \left \langle dP_{\eta}(s), \int_{-\frac{1}{2}}^{\frac{1}{2}} \frac{1}{a_{per} + s c_{per}} \right \rangle^{-1}.$$

The proof thus consists in inserting expansion (\ref{pushexpand}) in this explicit expression and identifying successively the first three dominant orders.\\

Using (\ref{pushexpand}), we write
\begin{eqnarray*}
(a_{\eta}^{*})^{-1} &=& \int_{-\frac{1}{2}}^{\frac{1}{2}} \frac{1}{a_{per}} + \eta \left \langle d\bar{P}_1(s), \int_{-\frac{1}{2}}^{\frac{1}{2}} \frac{1}{a_{per}+ sc_{per}} \right \rangle   + \eta^2 \left \langle d\bar{P}_2(s), \int_{-\frac{1}{2}}^{\frac{1}{2}} \frac{1}{a_{per}+ sc_{per}} \right \rangle  \\
&& + o(\eta^2) \\
&& \hspace{-2cm} = (a_{per}^{*})^{-1} \left(1+\eta a_{per}^{*} \left \langle d\bar{P}_1(s), \int_{-\frac{1}{2}}^{\frac{1}{2}} \frac{1}{a_{per}+ sc_{per}} \right \rangle + \eta^2 a_{per}^{*} \left \langle d\bar{P}_2(s), \int_{-\frac{1}{2}}^{\frac{1}{2}} \frac{1}{a_{per}+ sc_{per}} \right \rangle  \right) \\
&& + o(\eta^2).
\end{eqnarray*}
This yields the expansion
\begin{eqnarray} \label{devexactgen}
\begin{aligned}
a_{\eta}^{*} =& a_{per}^{*} - \eta (a_{per}^{*})^2\left \langle d\bar{P}_1(s), \int_{-\frac{1}{2}}^{\frac{1}{2}} \frac{1}{a_{per}+ sc_{per}} \right \rangle \\ 
& + \eta^2 (a_{per}^{*})^3\left \langle d\bar{P}_1(s), \int_{-\frac{1}{2}}^{\frac{1}{2}} \frac{1}{a_{per}+ sc_{per}} \right \rangle^2 \\
& - \eta^2 (a_{per}^{*})^2 \left \langle d\bar{P}_2(s), \int_{-\frac{1}{2}}^{\frac{1}{2}} \frac{1}{a_{per}+ sc_{per}} \right \rangle + o(\eta^2).
\end{aligned}
\end{eqnarray}

We now devote the rest of the proof to verifying that the coefficients of $\eta$ and $\eta^2$ in (\ref{devexactgen}) are indeed obtained as the limit as $N \rightarrow \infty$ of $a_1^{*,N}$ and $a_2^{*,N}$ defined generally by (\ref{order1}) and (\ref{order2}) respectively, in this particular one-dimensional setting.\\

The function $w^{1,s,0,N}$ generally defined by (\ref{sdefect}) satisfies here
\begin{eqnarray}\label{sdefect1d} 
\left  \{
\begin{aligned}
& -\frac{d}{dx}\left(a_1^{s,0}(\frac{d}{dx} w_i^{1,s,0,N} + 1 ) \right) = 0 \quad \, \mathrm{in}  \, ]-\frac{N}{2},\frac{N}{2}[,\\
& w_i^{1,s,0,N} \, N-\mathrm{periodic}.
\end{aligned}
\right.
\end{eqnarray} 

We easily compute using (\ref{sdefect1d}):
\begin{eqnarray*}
a_{1}^{s,0} (\frac{d}{dx} w^{1,s,0,N}+1) &=& N \left(\int_{-\frac{N}{2}}^{\frac{N}{2}} \frac{1}{a_{per}+s\mathds{1}_{[-\frac{1}{2},\frac{1}{2}]} c_{per}}\right)^{-1} \\
&=& N \left(N (a_{per}^*)^{-1}- f(s) \right)^{-1}\\
&=& a_{per}^{*}+ \frac{(a_{per}^{*})^2}{N}f(s) + \frac{(a_{per}^{*})^3}{N^2 } f(s)^2+ o(N^{-2}),
\end{eqnarray*}

where $\displaystyle f(s) = \int_{-\frac{1}{2}}^{\frac{1}{2}} \frac{sc_{per}}{a_{per} (a_{per}+sc_{per})}$.

Thus $a_1^{*,N}$ defined generally by (\ref{order1}) takes here the form
\begin{eqnarray*}
a_1^{*,N} &=& \left \langle d\bar{P}_1(s), \int_{-\frac{N}{2}}^{\frac{N}{2}} a_{1}^{s,0} (\frac{d}{dx} w^{1,s,0,N}+1) \right \rangle\\
&=&N a_{per}^{*} \left \langle d\bar{P}_1(s), 1 \right \rangle + (a_{per}^{*})^2 \left \langle d\bar{P}_1(s), f(s) \right \rangle +o(1).  
\end{eqnarray*}
We know from Lemma \ref{distreasy} that $\left \langle d\bar{P}_1(s), 1 \right \rangle = 0$, whence
\begin{eqnarray} \label{ainf1d}
a_1^{*,N} \underset{N \rightarrow \infty}{\rightarrow} \bar{a}_1^* = (a_{per}^{*})^2 \left \langle d\bar{P}_1(s), f(s) \right \rangle.
\end{eqnarray}

Likewise, we compute from (\ref{sdefect2}), for $k \in \llbracket -\frac{N-1}{2}, \frac{N-1}{2}\rrbracket$,
\begin{eqnarray*}
a_{2}^{s,t,0,k} (\frac{d}{dx} w^{2,s,t,0,k,N}+1) &=& N \left(\int_{-\frac{N}{2}}^{\frac{N}{2}} \frac{1}{a_{per}+s\mathds{1}_{[-\frac{1}{2},\frac{1}{2}]}c_{per}+t\mathds{1}_{[k-\frac{1}{2},k+\frac{1}{2}]}c_{per}}\right)^{-1} \\
&& \hspace{-3cm} =  N \left(N (a_{per}^*)^{-1}-\int_{-\frac{1}{2}}^{\frac{1}{2}} \frac{sc_{per}}{a_{per} (a_{per}+sc_{per})} - \int_{-\frac{1}{2}}^{\frac{1}{2}} \frac{tc_{per}}{a_{per} (a_{per}+tc_{per})} \right)^{-1}\\
&& \hspace{-3cm} = N \left(N (a_{per}^*)^{-1}-f(s)-f(t)\right)^{-1}.
\end{eqnarray*}
Then
$$a_{2}^{s,t,0,k} (\frac{d}{dx} w^{2,s,t,0,k,N}+1) =a_{per}^{*}+ \frac{(a_{per}^{*})^2}{N} (f(s)+f(t))+ \frac{(a_{per}^{*})^3}{N^2 }(f(s)+f(t))^2 + o(N^{-2}).$$

Notice that this expression is independent of $k$ (and so of the distance between the two defects), so that $a_2^{*,N}$ defined by (\ref{order2}) here reads
\begin{eqnarray} \label{interm}
\begin{aligned}
a_2^{*,N} =& \frac{N (N-1)}{2}  \left \langle d\bar{P}_1(s) d\bar{P}_1(t), a_{per}^{*}+ \frac{(a_{per}^{*})^2}{N} (f(s)+f(t))+ \frac{(a_{per}^{*})^3}{N^2 }(f(s)+f(t))^2 \right \rangle   \\
& + N \left \langle d\bar{P}_2(s), a_{per}^{*}+ \frac{(a_{per}^{*})^2}{N}f(s) \right \rangle  + o(1).
\end{aligned}
\end{eqnarray}

Since we know from Lemma \ref{distreasy} that $\left \langle d\bar{P}_1(s), 1 \right \rangle = 0$ and $\left \langle d\bar{P}_2(s), 1 \right \rangle = 0$, (\ref{interm}) reduces to
\begin{eqnarray*}
a_2^{*,N} =(a_{per}^{*})^3 \left \langle d\bar{P}_1(s), f(s) \right \rangle^2 + (a_{per}^{*})^2 \left \langle d\bar{P}_2(s), f(s) \right \rangle + o(1).
\end{eqnarray*}
Thus
\begin{eqnarray} \label{ainf2d}
a_2^{*,N}  \underset{N \rightarrow \infty}{\rightarrow} \bar{a}_2^* = (a_{per}^{*})^3 \left \langle d\bar{P}_1(s), f(s) \right \rangle^2 + (a_{per}^{*})^2 \left \langle d\bar{P}_2(s), f(s) \right \rangle.
\end{eqnarray}

Finally, since $\displaystyle f(s) = \int_{-\frac{1}{2}}^{\frac{1}{2}} \frac{1}{a_{per}} - \int_{-\frac{1}{2}}^{\frac{1}{2}} \frac{1}{a_{per} + s c_{per}}$, and $\left \langle d\bar{P}_1(s), 1 \right \rangle =  \left \langle d\bar{P}_2(s), 1 \right \rangle = 0$, we have
\begin{eqnarray} \label{simplify}
\left \langle d\bar{P}_1(s), f(s) \right \rangle = - \left \langle d\bar{P}_1(s), \int_{-\frac{1}{2}}^{\frac{1}{2}}\frac{1}{a_{per}+ sc_{per}} \right \rangle,
\end{eqnarray}
and
\begin{eqnarray} \label{simplify2}
\left \langle d\bar{P}_2(s), f(s) \right \rangle = - \left \langle d\bar{P}_2(s), \int_{-\frac{1}{2}}^{\frac{1}{2}}\frac{1}{a_{per}+ sc_{per}} \right \rangle.
\end{eqnarray}

In view of (\ref{devexactgen}), (\ref{ainf1d}), (\ref{ainf2d}), (\ref{simplify}) and (\ref{simplify2}), we have proved
$$ a_{\eta}^{*} = a_{per}^{*} + \eta \bar{a}_1^*  + \eta^2 \bar{a}_2^* + o(\eta^2).$$ 
\end{proof}

\subsection{A proof of the approach of Section \ref{heuristic} in a specific setting} \label{consistency}

The purpose of this final section is to prove that the formal approach of Section \ref{heuristic} is rigorous in a setting related to that of Corollary \ref{correl}.\\

More precisely, we assume that the random field $b_{\eta}$ satisfies the assumptions of Corollary \ref{correl}. These assumptions do not imply that the image measure $dP_{\eta}$ satisfies assumption (\ref{pushexpand}) which is at the heart of the approach of Section \ref{heuristic}, so that we have to impose that $dP_{\eta}$ additionally satisfies (\ref{pushexpand}). The following preliminary result then gives the necessary form of the expansion of the image measure $dP_{\eta}$.

\begin{lemma} \label{checkcor}
Assume that $b_{\eta}$ satisfies 
\begin{eqnarray}
b_{\eta}(x,\omega) = \sum_{k \in \mathbb{Z}^d} \mathds{1}_{Q+k}(x) B _{\eta}^k(\omega),
\end{eqnarray}
where the $B _{\eta}^k$ are i.i.d random variables, the distribution of which is given by a ``mother variable'' $B_{\eta}$ satisfying
\begin{eqnarray} 
&\forall \eta >0, \|B_{\eta}\|_{L^{\infty}(\Omega)} \leq M, \\
&B_{\eta} = \eta \bar{B}_0 + \eta^2 \bar{R}_0 + o(\eta^2) \quad \mathrm{weakly \; in\; } L^2(\Omega). \label{convforte}
\end{eqnarray} 

Assume further that the image measure $dP_{\eta}$ of $B_{\eta}$ satisfies (\ref{pushexpand}). Then
\begin{eqnarray} \label{expanspart2}
dP_{\eta} = \delta_0 - \eta \mathbb{E}(\bar{B}_0)\delta_0'+  \frac{\eta^2}{2} \mathbb{E}(\bar{B}_0^2)\delta_0'' - \eta^2 \mathbb{E}(\bar{R}_0) \delta_0'+ o(\eta^2) \; \mathrm{in} \; \mathcal{E}'(\mathbb{R}).
\end{eqnarray}
\end{lemma}

\begin{proof}
Firstly, notice that $\frac{B_{\eta}}{\eta}$ converges strongly to $\bar{B}_0$ in $L^2(\Omega)$ because of (\ref{convforte}). Now consider $\varphi \in \mathcal{D}(\RR)$. We have on the one hand
$$ \mathbb{E} \left(\frac{B_{\eta}^2}{\eta^2} \varphi(B_{\eta}) \right) \rightarrow  \mathbb{E} \left(\bar{B}_0^2\right) \varphi(0),$$
and on the other hand
$$ \mathbb{E}\left(B_{\eta}^2 \varphi(B_{\eta}) \right)= \eta \langle s^2 d\bar{P}_1 , \varphi \rangle + \eta^2 \langle s^2 d\bar{P}_2 , \varphi \rangle + o(\eta^2).$$

Thus $s^2 d\bar{P}_1 = 0$ and $s^2 d\bar{P}_2 = \mathbb{E} \left(\bar{B}_0^2\right)  \delta_0$ in $\mathcal{D}'(\RR)$. It is then well known that there exist $\gamma_1$, $\kappa_1$, $\gamma_2$, $\kappa_2$ in $\RR$ such that
$$d\bar{P}_1 = \gamma_1 \delta_0 + \kappa_1 \delta_0' \quad \mathrm{and} \quad d\bar{P}_2 = \gamma_2 \delta_0 + \kappa_2 \delta_0' + \frac{\mathbb{E} \left(\bar{B}_0^2\right)}{2} \delta_0''.$$ 

Lemma \ref{distreasy} implies $\gamma_1 = \gamma_2 = 0$. Then, we have
$$ \mathbb{E}(B_{\eta}) = \eta \mathbb{E}(\bar{B}_0) + \eta^2 \mathbb{E}(\bar{R}_0) + o(\eta^2)$$ and also
$$\mathbb{E}(B_{\eta}) = \eta \langle s d\bar{P}_1 , 1 \rangle +  \eta^2 \langle s d\bar{P}_2 , 1 \rangle + o(\eta^2).$$

Thus $\langle s d\bar{P}_1 , 1 \rangle =\mathbb{E}(\bar{B}_0) $ and $\langle s d\bar{P}_2 , 1 \rangle = \mathbb{E} \left(\bar{R}_0\right)$, from which we deduce $\kappa_1 = -\mathbb{E}(\bar{B}_0) $ and $\kappa_2 = -\mathbb{E}(\bar{R}_0)$. 

\end{proof}

Theorem \ref{rig} and Corollary \ref{correl} rigorously yield the second-order expansion 
$$ A_{\eta}^* = A_{per}^* + \eta \tilde{A_1^*} + \eta^2 \tilde{A_2^*} + o(\eta^2)$$
with $\tilde{A_1^*}$ and $\tilde{A_2^*}$ respectively defined by (\ref{tilda1}) and (\ref{tilda2}).

On the other hand, using (\ref{expanspart2}), Section \ref{heuristic} yields the formal expansion
$$ A_{\eta}^* = A_{per}^* + \eta \bar{A_1^*} + \eta^2 \bar{A_2^*} + o(\eta^2).$$
where $\bar{A_1^*}$ is the limit of the sequence $A_{1}^{*,N}$ defined by (\ref{order1}) or equivalently by (\ref{order1mod}), and $\bar{A_2^*}$ the limit (proved only up to extraction) of the sequence $A_{2}^{*,N}$ defined by (\ref{order2}).\\

The rest of this section is devoted to verifying that $\bar{A_1^*}$ coincides with $\tilde{A}_1^*$ and $\bar{A_2^*}$ coincides with $\tilde{A}_2^*$ in the specific setting of Lemma \ref{checkcor}.

\subsubsection{First-order term}

Using (\ref{expanspart2}), (\ref{order1mod}) reads  
\begin{eqnarray*} 
A_{1}^{*,N} e_i \cdot e_j = - \mathbb{E}(\bar{B}_0) \left \langle \delta_0'(s) , \int_{Q} sC_{per}(\nabla w_i^{1,s,0,N} + e_i) \cdot (e_j + \nabla \tilde{w}_{j}^0) \right \rangle,
\end{eqnarray*}
and we compute
\begin{eqnarray*} 
A_{1}^{*,N} = \mathbb{E}(\bar{B}_0) \int_{Q} C_{per}(\nabla w_i^{1,0,0,N} + e_i) \cdot (e_j + \nabla \tilde{w}_{j}^0).
\end{eqnarray*}
Setting $s=0$ in (\ref{sdefect}), it is clear that $w_i^{1,0,0,N}$ is equal to the periodic corrector $w_i^{0}$. Then
\begin{eqnarray} \label{order1consist} 
A_{1}^{*,N} = \mathbb{E}(\bar{B}_0) \int_{Q} C_{per}(\nabla w_i^{0} + e_i) \cdot (e_j + \nabla \tilde{w}_{j}^0). 
\end{eqnarray}

Clearly $A_{1}^{*,N}$ does not depend on $N$ and its limit is then
\begin{eqnarray} \label{order1consistlim} 
\bar{A_{1}^{*}} = \mathbb{E}(\bar{B}_0) \int_{Q} C_{per}(\nabla w_i^{0} + e_i) \cdot (e_j + \nabla \tilde{w}_{j}^0). 
\end{eqnarray}

We recognize in the right-hand side of (\ref{order1consistlim}) the first-order coefficient in (\ref{other}), which we know from Remark \ref{trick} is equivalent to (\ref{tilda1}). Theorem \ref{rig} therefore shows that the first-order expansion 
$$A_{\eta}^* = A_{per}^* + \eta \bar{A_{1}^{*}} + o(\eta)$$
is correct with the values of the coefficients given by our formal approach of Section \ref{heuristic}.\\

We now proceed similarly with the second-order coefficient.

\subsubsection{Second-order term}

Using the adjoint cell problems (\ref{celladj}) in (\ref{order2}) as in the proof of Proposition \ref{conva1}, let us first rewrite
\begin{eqnarray} \label{order2mod}
\begin{aligned}
A_{2}^{*,N} e_i \cdot e_j =& \sum_{k \in \mathcal{T}_N, k \neq 0} \left \langle d\bar{P}_1(s) d\bar{P}_1(t), \int_{Q} s C_{per}\nabla w_i^{2,s,t,0,k,N}\cdot (e_j + \nabla \tilde{w}_{j}^0) \right \rangle  \\
&+ \left \langle d\bar{P}_2(s), \int_{Q}  s C_{per} (\nabla w_i^{1,s,0,N} + e_i)\cdot (e_j + \nabla \tilde{w}_{j}^0) \right \rangle.
\end{aligned}
\end{eqnarray}
Inserting (\ref{expanspart2}) in (\ref{order2mod}), we start by focusing on
\begin{eqnarray*}
\left \langle d\bar{P}_2(s), \int_{Q}  s C_{per} (\nabla w_i^{1,s,0,N} + e_i )\cdot (e_j + \nabla \tilde{w}_{j}^0) \right \rangle = & & \\
& & \hspace{-9cm}  \frac{1}{2}\left(\mathbb{E}(\bar{B}_0) \right)^2 \left \langle \delta_0''(s), \int_{Q}  s C_{per} (\nabla w_i^{1,s,0,N} + e_i)\cdot (e_j + \nabla \tilde{w}_{j}^0) \right \rangle   \\
& & \hspace{-9cm} - \, \mathbb{E}(\bar{R}_0) \left \langle \delta_0'(s), \int_{Q}  s C_{per} (\nabla w_i^{1,s,0,N} + e_i)\cdot (e_j + \nabla \tilde{w}_{j}^0) \right \rangle.
\end{eqnarray*}
Denoting by $\partial_s w_i^{1,0,0,N}$, the first derivative of $w_i^{1,s,0,N}$ evaluated at $s=0$, we compute
\begin{eqnarray*}
\begin{aligned}
& \left \langle d\bar{P}_2(s), \int_{Q}  s C_{per} (\nabla w_i^{1,s,0,N} + e_i)\cdot (e_j + \nabla \tilde{w}_{j}^0) \right \rangle  & & \\
& = \mathbb{E}(\bar{B}_0^2) \int_{Q} C_{per} \nabla \partial_s w_i^{1,0,0,N}  \cdot (\nabla \tilde{w}_j^0 + e_j) + \mathbb{E}(\bar{R}_0) \int_{Q} C_{per} (\nabla  w_i^{1,0,0,N} + e_i) \cdot (\nabla \tilde{w}_j^0 + e_j)  \\
& = \mathbb{E}(\bar{B}_0^2) \int_{Q} C_{per} \nabla \partial_s w_i^{1,0,0,N}  \cdot (\nabla \tilde{w}_j^0 + e_j) + \mathbb{E}(\bar{R}_0) \int_{Q} C_{per} (\nabla  w_i^{0} + e_i) \cdot (\nabla \tilde{w}_j^0 + e_j).
\end{aligned}
\end{eqnarray*}

 It follows from (\ref{sdefect}) that $\partial_s w_i^{1,0,0,N}$ solves
\begin{eqnarray} \label{partialsw}
\left  \{
\begin{aligned}
& -\mathrm{div}(A_{per} \nabla \partial_s w_i^{1,0,0,N}) = \mathrm{div}\left(\mathds{1}_Q C_{per}(\nabla w_i^0 + e_i)\right) \quad \, \mathrm{in}  \, I_N,\\
& \partial_s w_i^{1,0,0,N} \, (N\ZZ)^d-\mathrm{periodic}. \\
\end{aligned}
\right.
\end{eqnarray}

Applying Lemma \ref{conv1} to (\ref{partialsw}), we deduce that $\nabla \partial_s w_i^{1,0,0,N}$ converges in $L^2(Q)$, when $N \rightarrow \infty$, to $\nabla t_i$ defined by (\ref{ti}) in Corollary \ref{correl}. Consequently, 
\begin{eqnarray} \label{convterm2}
\left \langle d\bar{P}_2(s), \int_{Q}  s C_{per} \left(\nabla w_i^{1,s,0,N} + e_i \right)\cdot (e_j + \nabla \tilde{w}_{j}^0) \right \rangle \underset{N \rightarrow \infty}{\rightarrow}&&  \nonumber \\
& & \hspace{-10cm}  \mathbb{E}(\bar{B}_0^2) \int_{Q} C_{per} \nabla t_i  \cdot (\nabla \tilde{w}_j^0 + e_j) + \mathbb{E}(\bar{R}_0) \int_{Q} C_{per} \left(\nabla  w_i^{0} + e_i \right) \cdot (\nabla \tilde{w}_j^0 + e_j). 
\end{eqnarray}

Next, we address
\begin{eqnarray*} 
\sum_{k \in \mathcal{T}_N, k \neq 0} \left \langle d\bar{P}_1(s) d\bar{P}_1(t), \int_{Q} s C_{per}\nabla w_i^{2,s,t,0,k,N}\cdot (e_j + \nabla \tilde{w}_{j}^0) \right \rangle &=&  \\
& &\hspace{-10cm}\mathbb{E}(\bar{B}_0)^2\sum_{k \in \mathcal{T}_N, k \neq 0} \left \langle  \delta_0'(s) \delta_0'(t), \int_{Q} s C_{per}\nabla w_i^{2,s,t,0,k,N}\cdot (e_j + \nabla \tilde{w}_{j}^0) \right \rangle.
\end{eqnarray*}

Denoting by $\partial_t w_i^{2,0,0,0,k,N}$ the first derivative of $w_i^{2,s,t,0,k,N}$ with respect to $t$ evaluated at $s=t=0$, we have
\begin{eqnarray} \label{termedur}
\sum_{k \in \mathcal{T}_N, k \neq 0} \left \langle d\bar{P}_1(s) d\bar{P}_1(t), \int_{Q} s C_{per}\nabla w_i^{2,s,t,0,k,N}\cdot (e_j + \nabla \tilde{w}_{j}^0) \right \rangle &=& \nonumber \\
& &\hspace{-10cm}\mathbb{E}(\bar{B}_0)^2\sum_{k \in \mathcal{T}_N, k \neq 0}  \int_{Q}  C_{per}\nabla \partial_t w_i^{2,0,0,0,k,N}\cdot (e_j + \nabla \tilde{w}_{j}^0).
\end{eqnarray}

It follows from (\ref{sdefect2}) that $\partial_t w_i^{2,0,0,0,k,N}$ solves
\begin{eqnarray} \label{partialsw2}
\left  \{
\begin{aligned}
& -\mathrm{div}(A_{per} \nabla \partial_t w_i^{2,0,0,0,k,N}) = \mathrm{div}\left(\mathds{1}_{Q+k} C_{per}(\nabla w_i^0 + e_i)\right) \quad \, \mathrm{in}  \, I_N,\\
& \partial_t w_i^{2,0,0,0,k,N} \, (N\ZZ)^d-\mathrm{periodic}. \\
\end{aligned}
\right.
\end{eqnarray}

Defining $d_i^N = \sum_{k \in \mathcal{T}_N} \partial_t w_i^{2,0,0,0,k,N}$, it is easy to see that $d_i^N$ is a $\ZZ^d$-periodic function that solves 
\begin{eqnarray} \label{din}
\left  \{
\begin{aligned}
& -\mathrm{div}(A_{per} \nabla d_i^N) = \mathrm{div}\left(C_{per}(\nabla w_i^0 + e_i)\right) \quad \, \mathrm{in}  \, Q, \\
& d_i^N \, \ZZ^d-\mathrm{periodic}.
\end{aligned}
\right.
\end{eqnarray}

Since problem (\ref{din}) has a unique solution up to an additive constant, $\nabla d_i^N = \nabla s_i$ where $s_i$ is defined by (\ref{si}) in Corollary \ref{correl}.\\

Finally, comparing (\ref{partialsw}) to (\ref{partialsw2}) for $k=0$, we find that $\nabla \partial_t w_i^{2,0,0,0,0,N}$ is equal to $\nabla \partial_s w_i^{1,0,0,N}$ and then also converges in $L^2(Q)$ to $\nabla t_i$ when $N \rightarrow \infty$.\\

Then, starting from (\ref{termedur}), 
\begin{eqnarray} \label{convterm1}
\sum_{k \in \mathcal{T}_N, k \neq 0} \left \langle d\bar{P}_1(s) d\bar{P}_1(t), \int_{Q} s C_{per}\nabla w_i^{2,s,t,0,k,N}\cdot (e_j + \nabla \tilde{w}_{j}^0) \right \rangle && \nonumber \\
& & \hspace{-11cm}= \left(\mathbb{E}(\bar{B}_0)\right)^2 \int_{Q} C_{per} \sum_{k \in \mathcal{T}_N}  \nabla \partial_t w_i^{2,0,0,0,k,N}\cdot (e_j + \nabla \tilde{w}_{j}^0) \nonumber \\
& & \hspace{-10cm} - \left(\mathbb{E}(\bar{B}_0)\right)^2 \int_{Q}C_{per} \nabla \partial_t w_i^{2,0,0,0,0,N} \cdot (e_j + \nabla \tilde{w}_{j}^0) \nonumber \\
& & \hspace{-11cm} = \left(\mathbb{E}(\bar{B}_0)\right)^2 \int_{Q} C_{per} \nabla s_i \cdot (e_j + \nabla \tilde{w}_{j}^0) - \left(\mathbb{E}(\bar{B}_0)\right)^2 \int_{Q} C_{per} \nabla \partial_t w_i^{2,0,0,0,0,N} \cdot (e_j + \nabla \tilde{w}_{j}^0) \nonumber \\
& & \hspace{-11cm} \underset{N \rightarrow \infty}{\rightarrow} \left(\mathbb{E}(\bar{B}_0)\right)^2 \int_{Q} C_{per} \nabla s_i \cdot (e_j + \nabla \tilde{w}_{j}^0) - \left(\mathbb{E}(\bar{B}_0)\right)^2 \int_{Q} C_{per} \nabla t_i \cdot (e_j + \nabla \tilde{w}_{j}^0)
\end{eqnarray}

It entails from (\ref{order2mod}), (\ref{convterm2}) and (\ref{convterm1}) that $A_2^{*,N}$ converges to a limit $\bar{A_2^*}$ defined by
\begin{eqnarray*}
\bar{A_2^*} e_i \cdot e_j &=& \mathbb{E}(\bar{R}_0) \int_{Q} C_{per} \left(\nabla  w_i^{0} + e_i \right) \cdot (\nabla \tilde{w}_j^0 + e_j) + Var(\bar{B}_0)\int_{Q} C_{per} \nabla t_i \cdot (\nabla \tilde{w}_j^0 + e_j)\\
& & \hspace{5cm} \quad \quad \quad + \left(\mathbb{E}(\bar{B}_0)\right)^2 \int_{Q} \nabla s_i \cdot (e_j + \nabla \tilde{w}_{j}^0).
\end{eqnarray*}
$\bar{A_2^*}$ is equal to the second-order term given by (\ref{expans2correl}) in Corollary \ref{correl} since we deal with independent random variables in each cell of $\ZZ^d$. Thus the second-order expansion
$$A_{\eta}^* = A_{per}^* + \eta  \bar{A_1^*} + \eta^2  \bar{A_2^*} + o(\eta^2)$$
derived from the formal approach of Section \ref{heuristic} is correct in this specific setting.

\end{document}